\documentclass[11pt]{article}

\RequirePackage[amsthm,amsmath]{imsart}

 \parskip        \baselineskip
\usepackage[mathscr]{eucal}
\usepackage{graphicx}
\usepackage{latexsym}
\usepackage{amssymb}
\usepackage{psfrag}
\usepackage{epsfig}
\DeclareMathAlphabet{\mathpzc}{OT1}{pzc}{m}{it}

\include{psbox}

\usepackage{amsfonts}

\usepackage{graphicx}

\usepackage{amsfonts}

\usepackage{color}

\begin{document}
\bibliographystyle{plain}

\newtheorem{Theorem}{\bf Theorem}[section]
\newtheorem{lemma}[Theorem]{\bf Lemma}
\newtheorem{example}{\bf Example}[section]
\newtheorem{definition}{\bf Definition}[section]

\newtheorem{Corollary}[Theorem]{\bf Corollary}
\newtheorem{remark}{\bf Remark}[section]

\newtheorem{Assumption}{Assumption}
\newtheorem{condition}{\bf Condition}[section]
\newtheorem{Proposition}[Theorem]{\bf Proposition}
\newtheorem{definitions}{\bf Definition}[section]
\numberwithin{equation}{section}
\newcommand{\skp}{\vspace{\baselineskip}}
\newcommand{\noi}{\noindent}
\newcommand{\ti}{\tilde}
\newcommand{\osc}{\mbox{osc}}

\newcommand{\eps}{\varepsilon}
\newcommand{\e}{\varepsilon}
\newcommand{\del}{\delta}
\newcommand{\Rk}{\mathbb{R}^k}
\newcommand{\R}{\mathbb{R}}
\newcommand{\spa}{\vspace{.2in}}
\newcommand{\V}{\mathcal{V}}
\newcommand{\E}{\mathbb{E}}
\newcommand{\I}{\mathbb{I}}
\newcommand{\p}{\mathbb{P}}
\newcommand{\PP}{\mathcal{P}}
\newcommand{\vs}{\varsigma}

\newcommand{\QQ}{\mathbb{Q}}

\newcommand{\lan}{\langle}
\newcommand{\ran}{\rangle}
\def\wh{\widehat}
\newcommand{\defn}{\stackrel{def}{=}}
\newcommand{\txb}{\tau^{\epsilon,x}_{B^c}}
\newcommand{\tyb}{\tau^{\epsilon,y}_{B^c}}
\newcommand{\tilxb}{\tilde{\tau}^\eps_1}
\newcommand{\pxeps}{\mathbb{P}_x^{\eps}}
\newcommand{\non}{\nonumber}
\newcommand{\dist}{\mbox{dist}}

\newcommand{\tilyb}{\tilde{\tau}^\eps_2}
\newcommand{\beq}{\begin{eqnarray*}}
\newcommand{\eeq}{\end{eqnarray*}}
\newcommand{\beqn}{\begin{eqnarray}}
\newcommand{\eeqn}{\end{eqnarray}}
\newcommand{\ink}{\rule{.5\baselineskip}{.55\baselineskip}}

\newcommand{\bt}{\begin{theorem}}
\newcommand{\et}{\end{theorem}}
\newcommand{\deps}{\Delta_{\eps}}

\newcommand{\be}{\begin{equation}}
\newcommand{\ee}{\end{equation}}
\newcommand{\ac}{\mbox{AC}}
\newcommand{\BB}{\mathbb{B}}
\newcommand{\VV}{\mathbb{V}}
\newcommand{\DD}{\mathbb{D}}
\newcommand{\KK}{\mathbb{K}}
\newcommand{\HH}{\mathbb{H}}
\newcommand{\TT}{\mathbb{T}}
\newcommand{\CC}{\mathbb{C}}
\newcommand{\SSS}{\mathbb{S}}
\newcommand{\EE}{\mathbb{E}}
\newcommand{\Lip}{\mbox{Lip}}

\newcommand{\RR}{{\mathbb{R}}}

\newcommand{\clg}{\mathcal{G}}
\newcommand{\clb}{\mathcal{B}}
\newcommand{\cls}{\mathcal{S}}
\newcommand{\clc}{\mathcal{C}}
\newcommand{\cle}{\mathcal{E}}
\newcommand{\clv}{\mathcal{V}}
\newcommand{\clu}{\mathcal{U}}
\newcommand{\clr}{\mathcal{R}}
\newcommand{\clt}{\mathcal{T}}
\newcommand{\cll}{\mathcal{L}}

\newcommand{\cli}{\mathcal{I}}
\newcommand{\clp}{\mathcal{P}}
\newcommand{\cla}{\mathcal{A}}
\newcommand{\clf}{\mathcal{F}}
\newcommand{\clh}{\mathcal{H}}
\newcommand{\N}{\mathbb{N}}
\newcommand{\Q}{\mathbb{Q}}

\newcommand{\curvz}{{\bf \mathpzc{z}}}
\newcommand{\curvx}{{\bf \mathpzc{x}}}
\newcommand{\curvi}{{\bf \mathpzc{i}}}
\newcommand{\curvs}{{\bf \mathpzc{s}}}

\newcommand{\tac}{\mbox{\scriptsize{AC}}}
\newcommand{\beginsec}{
\setcounter{lemma}{0} \setcounter{theorem}{0}
\setcounter{corollary}{0} \setcounter{definition}{0}
\setcounter{example}{0} \setcounter{proposition}{0}
\setcounter{condition}{0} \setcounter{assumption}{0}
\setcounter{remark}{0} }

\numberwithin{equation}{section} 

\begin{frontmatter}
\title{ Long Time Results for a Weakly Interacting Particle System in Discrete Time.}

 \runtitle{Weakly interacting particles.}

\begin{aug}
\author{Amarjit Budhiraja and Abhishek Pal Majumder\\ \ \\
}
\end{aug}

\today

\skp

\begin{abstract}
\noi
We study long time behavior of a discrete time weakly interacting particle system, and the corresponding nonlinear Markov process in $\R^d$, described in terms of a general stochastic evolution equation. In a setting where the state space of the particles is compact such questions have been studied in previous works, however  for the case of an unbounded state space 
very few results are available.
Under suitable assumptions on the problem data we study several  time asymptotic properties of the $N$-particle system and the associated nonlinear Markov chain.  In particular we show that the evolution equation for the law of the nonlinear Markov chain   has a unique fixed point and starting from an arbitrary initial condition convergence to the fixed point occurs at an exponential rate.  The empirical measure $\mu_n^N$ of the $N$-particles at time $n$  is shown to converge  to the law $\mu_n$ of the nonlinear Markov process at time $n$, in the Wasserstein-1 distance, in $L^1$,
as $N \to \infty$, {\em uniformly} in $n$. Several consequences of this uniform convergence are studied, including the interchangeability of the limits $n \to \infty$ and $N \to \infty$ and  the propagation of chaos property at $n=\infty$. Rate of convergence 
of $\mu_n^N$ to $\mu_n$ is studied
by establishing uniform in time polynomial and exponential probability concentration estimates.

\noi {\bf AMS 2010 subject classifications:} Primary 60J05, 60K35, 60F10.

\noi {\bf Keywords:} Weakly interacting particle system, propagation of chaos,  nonlinear Markov chains, Wasserstein distance, McKean-Vlasov equations, exponential concentration estimates, transportation inequalities, metric entropy, stochastic difference equations, long time behavior, uniform concentration estimates.
\end{abstract}

\end{frontmatter}

\section{Introduction}
\label{sec:intro}

Stochastic dynamical systems that model the evolution of a large collection of weakly interacting particles have long been studied in statistical mechanics (cf. \cite{sznitman1991topics, kolk}
and references therein).  In recent years such models have been considered in many other application areas as well, some examples include,
chemical and biological systems( e.g. biological aggregation, chemotactic response dynamics\cite{schweitzer2007brownian, Ste, friedman2009asymptotic}),
mathematical finance (e.g. mean field games\cite{GLL,CaDe}, default clustering in large portfolios\cite{GSS}), 
social sciences (e.g. opinion dynamics models \cite{como2009scaling, gomez2012bounded}), communication systems (\cite{GHK,AFRT,graham2009interacting}) etc.
Starting from the work of Sznitman\cite{sznitman1991topics} there has been an extensive body of work that studies law of large number behavior (Propagation of Chaos),
central limit theory (normal fluctuations from the mean) and large deviation principles for such models.
All of these results concern the behavior of the system over a finite time horizon.  In many applications the time asymptotic behavior of the system is of central concern.  
For example, stability of a communication system, steady state aggregation and self organization in biological and chemical systems,  long term consensus formation
mechanisms in opinion dynamics modeling, particle based approximation methods for invariant measures all rely on a careful analysis of the time asymptotic behavior of such systems.
Such behavior for special families of weakly interacting particle systems has been considered by several authors.  In \cite{del2011concentration} the authors give
general sufficient conditions for a family of discrete time systems for uniform in time exponential probability concentration estimates to hold. These conditions formulated in terms of Dobrushin's coefficient  are not very restrictive when the state space of the particles is compact, however they are hard to verify for settings with an unbounded state space.
In \cite{budhiraja2011discrete} a discrete time model with a compact state space for chemotactic cell response dynamics was studied.  Several time asymptotic results, including uniform in time law of large numbers, exponential stability of the associated nonlinear Markov process and uniform in time convergence of a particle based simulation scheme, were established.
For the setting of an unbounded state space and in continuous time, there have been several recent interesting works on  granular media equations \cite{malrieu2003convergence, bolley2007quantitative, carrillo2003kinetic} which establish uniform in time propagation of chaos, time uniform convergence of simulation schemes and uniform in time exponential concentration estimates.  

In the current work we consider a discrete time weakly interacting particle system and the corresponding nonlinear Markov process in $\R^d$, described in terms of a general stochastic evolution equation. Denoting by $X_{n}^{i}\equiv X_{n}^{i,N}$ the state of the $i$-th particle $(i=1,...,N)$ at time instant $n$, the evolution is given as 
\begin{eqnarray}
X_{n+1}^{i} &=& AX_{n}^{i} + \delta f(X_{n}^{i}, \mu_{n}^{N},\epsilon_{n+1}^{i}),\text{\hspace{10 mm}} \quad i=1,..,N, \quad n\in\mathbb{N}_{0}\label{nlps}
\end{eqnarray}
Here $\mu_{n}^{N}:= \frac{1}{N}\sum_{i=1}^N \delta_{X_{n}^{i}}$ is the empirical measure of the particle values at time instant $n$, $A$ is a $d\times d$ matrix,  $\delta $ is a small parameter, $\{\epsilon_{n}^{i},i=1,...,N, \quad n\geq 1 \}$ is an i.i.d array of $\mathbb{R}^{m}$ valued random variables with common probability law $\theta$ and $f :\mathbb{R}^{d}\times\mathcal{P}(\mathbb{R}^{d})\times\mathbb{R}^{m} \to \mathbb{R}^{d}$ is a measurable function. Also, $\{X_{0}^{i},i=1,...,N\}$ are taken to be exchangeable with common distribution $\mu_{0}$.
As will be seen in Section \ref{sec:mainresults}, the following   nonlinear Markov chain will correspond to the $N\to\infty$ limit of  \eqref{nlps}.
\begin{eqnarray}
X_{n+1} = AX_{n} + \delta f(X_{n},\mu_{n},\epsilon_{n+1}), \quad  \mathcal{L}(X_{n})= \mu_{n},\;  n \in \mathbb{N}_0. \label{nls}
\end{eqnarray}
where throughout we denote by $\mathcal{L}(X)$ the probability distribution of a  random variable $X$ with values in some Polish space $S$.
Under conditions on $f,\theta,\delta$  and $A$ we study several long time properties of the $N$-particle system and the associated nonlinear Markov chain.
Our starting point is the evolution equation for the law of the nonlinear Markov chain given in \eqref{recursion}.  We show in Theorem \ref{th2} that under conditions, that include a Lipschitz property of $f$
with the Wasserstein-1($W1$) distance on the space of probability measures (Assumptions \ref{assumA} and \ref{assumB}), contractivity of $A$ (Assumption \ref{assumC}) and $\delta$ being sufficiently small, \eqref{recursion} has a unique fixed point and starting from an arbitrary initial condition convergence to the fixed point occurs at an exponential rate.  Using this result we next argue in Theorem \ref{th3} that under an additional integrability condition (Assumption \ref{assumD}), as $N \to \infty$, the empirical measure $\mu_n^N$ of the $N$-particles at time $n$ converges to the law $\mu_n$ of the nonlinear Markov process at time $n$, in the W1 distance, in $L^1$,
 {\em uniformly} in $n$.  This result in particular shows that the $W1$ distance between $\mu_n^N$ and the unique fixed point $\mu_{\infty}$ of \eqref{recursion} converges to zero as $n \to \infty$ and $N\to \infty$ in any order.  This result is key in developing particle based numerical schemes for approximating the fixed point of the evolution
equation \eqref{recursion}.  We next show that under an irreducibility condition on the underlying Markovian dynamics (Assumption \ref{assumE}) the unique invariant measure $\Pi^N_{\infty}$ of the $N$-particle dynamics is $\mu_{\infty}$-chaotic, namely as $N \to \infty$, the projection of $\Pi^N_{\infty}$ on the first $k$-coordinates converges to $\mu_{\infty}^{\otimes k}$ for every $k \ge 1$.  This propagation of chaos property all the way to $n=\infty$ crucially relies on the uniform in time convergence of $\mu_n^N$ to $\mu_{\infty}$.
The next three results study the rate of this uniform convergence by developing suitable probability concentration estimates.  The first result (Theorem \ref{th5}), under an assumption of polynomial moments 
on the initial data and noise sequence (Assumption \ref{assumD}) establishes a corresponding uniform in time polynomial concentration bound.  The proof relies on an idea of restricting measures  to a compact set and estimates on metric entropy introduced in \cite{bolley2007quantitative} (see also \cite{villani2003topics}).  The basic idea is to first obtain a concentration bound for the   $\mathcal{W}_{1}$ distance between the truncated law and its corresponding empirical law in a compact ball of radius $R$ along with an estimate on the contribution from the region outside the ball and finally  optimize suitably over $R$.
The last two results are concerned with exponential concentration.  These impose much stronger integrability conditions on the problem data (Assumption \ref{assumG}).  The first considers the setting where the initial random variables form a general
exchangeable sequence and gives a concentration bound with an exponential decay rate of $N^{\frac{1}{d+2}}$.  The second result uses exponential concentration estimates for empirical measures of i.i.d. sequences based
on transportation inequalities from \cite{boissard2011simple, bolley2005weighted} (see also \cite{djellout2004transportation, gozlan2010transport, gozlan2007large, bolley2005weighted, bolley2007quantitative}) and considers the setting where the initial data is i.i.d.
In this case the concentration bound gives an exponential decay rate of order $N$.

The following notation will be used in this work. $\mathbb{R}^{d}$ will denote the $d$ dimensional Euclidean space with the usual Euclidean norm $|.|$. The set of natural numbers (resp. whole numbers) is denoted by $\mathbb{N}$ (resp. $\mathbb{N}_{0}$). Cardinality of a finite set $S$ is denoted by $|S|$. For a measurable space S, $\mathcal{P}(S)$ denotes the space of all probability measures on $S$.   For $x\in \mathbb{R}^{d}$, $\delta_{x}$ is the Dirac delta measure on $\mathbb{R}^{d}$ that puts a unit mass at location $x$. The space of real valued bounded measurable functions on $S$ is denoted as $BM(S)$. Borel $\sigma$ field on a metric space will be denoted as $\mathcal{B}(S)$. $\mathcal{C}_{b}(S)$ denotes the space of all bounded and continuous functions $f: S \to \mathbb{R}$. The supremum norm of a function $f:S\to \mathbb{R}$ is $||f||_{\infty}=\sup_{x\in S}|f(x)|$. When $S$ is a metric space, the Lipschitz seminorm of $f$ is defined by $||f||_{1}=\sup_{x\not= y}\frac{|f(x)-f(y)|}{d(x,y)}$ where $d$ is the metric on the space $S$. For a bounded Lipschitz function $f$ on $S$ we define $||f||_{BL}:=||f||_{1}+||f||_{\infty}$.
$\mbox{Lip}_{1}(S)$ (resp. $BL_{1}(S)$ )  denotes the class of Lipschitz (resp. bounded Lipschitz) functions $f :S\to \mathbb{R}$ with $||f||_{1}$ (resp. $||f||_{BL}$)
bounded by 1. Occasionally we will suppress $S$ from the notation and write $\mbox{Lip}_{1}$ and  $BL_{1}$ when clear from the context.  For a Polish space $S$, $\mathcal{P}(S)$ is equipped with the topology of weak convergence. A convenient metric metrizing this topology on $\mathcal{P}(S)$ is given as $\beta(\mu,\gamma) = \sup \{|\int fd\mu - \int fd\gamma|:||f||_{BL_{1}}\leq 1 \}$ for $ \mu,\gamma\in  \mathcal{P}(S)$. For a signed measure $\gamma$ on $\mathbb{R}^{d}$, we define $\langle f,\gamma\rangle:= \int f d\gamma$ whenever the integral makes sense. Let $\mathcal{P}_{1}(\mathbb{R}^{d})$ be the space of $\mu \in \mathcal{P}(\mathbb{R}^{d})$ such that $$||\mu||_{1}:= \int|x|d\mu(x) <\infty.$$ The space $\mathcal{P}_{1}(\mathbb{R}^{d})$ will be equipped with the Wasserstein-1 distance that is defined as follows:
$$\mathcal{W}_{1}(\mu_{0},\gamma_{0})  := \inf_{X,Y}E|X-Y|$$ where the infimum is taken over all  $\mathbb{R}^{d}$ valued random variables X,Y defined on a common probability space and where the marginals of X, Y are respectively $\mu_{0}$ and $\gamma_{0}$.
From Kantorovich-Rubenstein duality (cf. \cite{villani2003topics}) one sees the Wasserstein-1 is same as 
\begin{eqnarray}
\mathcal{W}_{1}(\mu_{0},\gamma_{0}) =\sup_{f\in \mbox{Lip}_{1}(\mathbb{R}^{d})}|\langle f,\mu_{0} - \gamma_{0}\rangle|.
\end{eqnarray}
For a signed measure $\mu$ on $(S,\mathcal{B}(S))$, the total variation norm of $\mu$ is defined as $|\mu|_{TV}:=\sup_{||f||_{\infty} \leq 1}\langle f,\mu \rangle$. Convergence in distribution of a sequence $\{X_{n}\}_{n\geq 1}$ of $S$ valued random variable to $X$ will be written as $X_{n}\Rightarrow X$.

A finite collection $\{Y_{1},Y_{2},..,Y_{N}\}$ of $S$ valued random variables is called exchangeable if $$\mathcal{L}(Y_{1},Y_{2},..,Y_{N}) = \mathcal{L}(Y_{\pi(1)},Y_{\pi(2)},..,Y_{\pi(N)})$$ for every permutation $\pi$ on the $N$ symbols $\{1,2,...,N\}$.  Let $\{Y_{i}^{N},i=1,..,N\}_{N\geq 1}$ be a collection of $S$ valued random variables, such that for every $N$, $\{Y_{1}^{N},Y_{2}^{N},..,Y_{N}^{N}\}$ is exchangeable. Let $\nu_{N} = \mathcal{L}(Y_{1}^{N},Y_{2}^{N},..,Y_{N}^{N})$. The sequence $\{\nu_{N}\}_{N\geq 1}$ is called $\nu$ -chaotic (cf.  \cite{sznitman1991topics})  for a $\nu \in \mathcal{P}(\mathcal{S})$, if for any $k\geq 1$, $f_{1},f_{2},...,f_{k} \in \mathcal{C}_{b}(\mathcal{S}),$ one has
\begin{eqnarray}
\lim_{N\to \infty} \langle f_{1}\otimes f_{2}\otimes ...\otimes f_{k}\otimes 1... \otimes1,\nu_{N} \rangle = \prod_{i=1}^{k}\langle f_{i} ,\nu \rangle. \label{chaotic}
\end{eqnarray}

Denoting the marginal distribution on first $k$ coordinates of $\nu_{N}$ by $\nu_{N}^{k}$, equation (\ref{chaotic}) says that, for every $k\geq 1,$ $\nu_{N}^{k} \rightarrow \nu^{\otimes k}$.

\section{Model description}
\label{sec:mod-desc}

Recall the system of $N$ interacting particles in $\mathbb{R}^{d}$ introduced in \eqref{nlps}.   Throughout we will assume that  $\{X_{0}^{i},i=1,...,N\}$ is exchangeable with common distribution $\mu_{0}$ where $\mu_{0} \in \PP_{1}(\RR^{d})$. Assumptions on $f,\theta,\delta$  and $A$ will be introduced shortly. Note that in the notation we have suppressed the dependence of the sequence $\{X_{n}^{i}\}$ on $N$.
Given $\rho \in \PP(\RR^{d})$ define a transition probability kernel $P^{\rho}: \RR^{d} \times \mathcal{B}(\RR^{d})\to [0,1]$ as  $$P^{\rho}(x,C) = \int_{\RR^{m}} 1_{[Ax+\delta f(x,\rho,z)\in C]} \theta(dz),\quad  \quad (x,C)\in \RR^{d} \times \mathcal{B}(\RR^{d}).$$
With an abuse of notation we will also denote by $P^{\rho}$ the map from $BM(\RR^{d})$ to itself, defined as 
$$P^{\rho} \phi(x) = \int_{\RR^{d}} \phi(y) P^{\rho}(x,dy), \quad \phi \in BM(\RR^{d}), x\in \RR^{d}.$$
For $\mu \in \PP(\RR^{d})$, let $\mu P^{\rho}\in \PP(\RR^{d})$ be defined as
$$\mu P^{\rho}(A) = \int_{\RR^{d}} P^{\rho}(x,A) \mu(dx), \quad A\in \mathcal{B}(\mathbb{R}^{d}).$$
Note that $\mu P^{\rho} = \mathcal{L}(AX + \delta f(X,\rho,\epsilon))$ when $\mathcal{L}(X,\epsilon) = \mu \otimes \theta$.

Under Assumptions \ref{assumA} and \ref{assumB} introduced in the next section it will follow that, for $\rho,\mu \in \mathcal{P}_{1}(\mathbb{R}^{d}),$ $\mu P^{\rho} \in \mathcal{P}_{1}(\mathbb{R}^{d})$ as well. Under these conditions, one can define $\Psi : \mathcal{P}_{1}(\mathbb{R}^{d}) \to \mathcal{P}_{1}(\mathbb{R}^{d})$ as 
\begin{equation}
\Psi (\mu) = \mu P^{\mu}.
\label{eq:544}
\end{equation}
Then the evolution of the law of the nonlinear Markov chain given in \eqref{nls} is given by the equation 
\begin{eqnarray}
\mu_{n+1}=\Psi(\mu_{n}).\label{recursion}
\end{eqnarray}
Using the above notation we see that $((X_{n}^{1},...,X_{n}^{N}),\mu_{n}^{N})$  is a $\mathbb(\mathbb{R}^{d})^{N} \times\mathcal{P}_{1}(\mathbb{R}^{d})$ valued discrete time Markov chain  defined recursively as follows. Let $X_{k}(N) \equiv (X_{k}^{1},X_{k}^{2},...,X_{k}^{N})$ and let $\mathcal{F}_{0} = \sigma\{X_{0}(N)\}.$ Then, for $k\geq 1$
\begin{eqnarray}
\begin{cases}
P(X_{k}(N) \in C |\mathcal{F}_{k-1}^N) = \bigotimes_{j=1}^{N}(\delta_{X_{k-1}^{j}} P^{\mu_{k-1}^{N}})(C)\hspace{3mm} \forall C \in \mathcal{B}(\mathbb{R}^{d})^{N} \nonumber \\
\mu_{k}^{N}= \frac{1}{N}\sum_{i=1}^N \delta_{X_{k}^{i}}\nonumber \\
\mathcal{F}_{k}^{N} =\sigma\{X_{k}(N)\} \vee \mathcal{F}_{k-1}^{N}. 
\end{cases}
\end{eqnarray}
\section{Main Results}\label{sec:mainresults}
Recall that $\{X_{0}^{i},i=1,...,N\}$ is assumed to be exchangeable with common distribution $\mu_{0}$ where $\mu_{0} \in \PP_{1}(\RR^{d})$.
We now introduce our assumptions on the nonlinearity.
\begin{Assumption}\label{assumA}
 $\int D(z) \theta(dz) = \sigma < \infty$, where
 $$  D(z) :=\sup_{x_{1} \neq x_{2}, \mu_{1} \neq \mu_{2}, x_1, x_2 \in \RR^d, \mu_1, \mu_2 \in \clp_1(\RR^d)}\frac{|f(x_{1},\mu_{1},z)-f(x_{2},\mu_{2},z)|}{|x_{1}-x_{2}|+\mathcal{W}_{1}(\mu_{1},\mu_{2})}, \; z \in \R^m .$$
\end{Assumption}
Note that the Assumption \ref{assumA} implies that  
\begin{eqnarray}
\sup_{(x,\mu) \in \R^d \times \clp_1(\R^d)}|f(x,\mu,z)|\leq (|x| +||\mu||_{1})D(z) + D_{1}(z), \; z \in \R^m\label{fbound}
\end{eqnarray}
 where $D_{1}(z) := |f(0,\delta_{0},z)|$. We impose the following condition on $D_{1}.$
\begin{Assumption}\label{assumB}
$\int D_{1}(z)\theta(dz) =c_{0} <\infty $.
\end{Assumption}
Our first result is a law of large numbers for $\mu_{n}^{N}$ as $N\to \infty$. Note that under Assumptions \ref{assumA} and \ref{assumB},  $\mu_{n}\in \mathcal{P}_{1}(\mathbb{R}^{d})$
for all $n\in \mathbb{N}_{0}$.
\begin{Theorem}\label{th1}
	Suppose Assumptions \ref{assumA} and \ref{assumB} hold and
 suppose that $E \mathcal{W}_{1}(\mu_{0}^{N},\mu_{0})\rightarrow 0$ as $N\to \infty$. Then, as $N\to \infty$, 
\begin{eqnarray}
 E \mathcal{W}_{1}(\mu_{n}^{N},\mu_{n}) \to 0 \label{pointwise}
\end{eqnarray} for all $n\geq 0$.
\end{Theorem}
\begin{remark}\label{rem1} Note that Theorem \ref{th1}  says that for all $n \ge 0$
$$ \lim_{N\to \infty}E\sup_{f\in \mbox{Lip}_{1}}|\langle f,\mu_{n}^{N} - \mu_{n}\rangle| = 0,$$ which in particular implies that $\mu_{n}^{N} \rightarrow \mu_{n}$ in probability, in $\mathcal{P}(\mathbb{R}^{d})$ (with the topology of weak convergence) as $N\to \infty$.
\end{remark}
Next we state a ``propagation of chaos" result which is an immediate consequence of Remark \ref{rem1} and exchangeability of of $\{X_{n}^{i}\}_{i=1}^{N}$.
\begin{Corollary}
Suppose Assumptions \ref{assumA} and \ref{assumB} hold. Then for any $k\geq 1$ and $n\in \mathbb{N}_{0},$  $\mathcal{L}(X_{n}^{1},X_{n}^{2},...,X_{n}^{k})\longrightarrow(\mathcal{L}(X_{n}) )^{\bigotimes k}$  as  $N\to \infty$. 
\end{Corollary}
For a $d\times d$ matrix B we denote its norm by $||B||,$ i.e. $||B|| = \sup_{x \in \mathbb{R}^{d}\setminus\{0\}} \frac{|Bx|}{|x|}$.
\begin{Assumption}\label{assumC}
$||A|| \leq e^{-\omega}$ for some $\omega > 0$.
\end{Assumption}
A measure $\mu^{*}\in \mathcal{P}_{1}(\mathbb{R}^{d})$ is called a fixed point for the evolution equation in (\ref{nls}), if $\mu^{*} = \Psi(\mu^{*}).$
Let $a_0 = \frac{1-e^{-\omega}}{2\sigma}$.
\begin{Theorem}\label{th2}
Suppose Assumptions \ref{assumA},\ref{assumB} and \ref{assumC} hold and that $\delta \in (0, a_0).$ Then there exists a unique fixed point $\mu_{\infty}$ of equation (\ref{recursion}). Furthermore, denoting for $\gamma \in \mathcal{P}_{1}(\mathbb{R}^{d}),$ 
$\mu_{n}[\gamma] = \underbrace{\Psi \circ \Psi ... \circ \Psi }_{\mbox{n \text{times}}}(\gamma)$, we have $$\limsup_{n \to \infty} \frac{1}{n} \log \mathcal{W}_{1}(\mu_{n}[\gamma],\mu_{\infty} ) < 0,$$ 
namely $\mu_{n}[\gamma]$ converges to $\mu_{\infty}$ as $n\to \infty$, at an exponential rate.
\end{Theorem}
Next, we study uniform in time (i.e $n$) convergence of $\mu_{n}^{N}$ to $\mu_{n}$  as the number of particles $N\to \infty$. For this, we will strengthen Assumptions \ref{assumA} and \ref{assumB} as follows.  
\begin{Assumption}\label{assumD}
 For some $\alpha>0$  $$E|X_{0}^{i}|^{1+\alpha} <\infty,\quad \int D(z)^{1+\alpha} \theta(dz) = \sigma_{1}(\alpha) <\infty,\quad\text{and}\quad \int (D_{1}(z))^{1+\alpha} \theta(dz) =c_{1}(\alpha) <\infty.$$
\end{Assumption}
\begin{Theorem}\label{th3} Suppose that Assumptions \ref{assumC} and \ref{assumD} hold. Also suppose that  $\delta \in (0,a_0).$ Then 
\begin{enumerate}
\item Given $\varepsilon >0,$ there exist $N_{0}(\varepsilon), n_{0}(\varepsilon) \in \mathbb{N}$ such that
$$E\mathcal{W}_{1}(\mu_{n}^{N},\mu_{n}) < \varepsilon\quad \quad \text{whenever}\quad n\geq n_{0}(\varepsilon),N\geq N_{0}(\varepsilon).$$
\item Suppose $E \mathcal{W}_{1}(\mu_{0}^{N},\mu_{0})\rightarrow 0$ as $N\to \infty$. Then 
$ \sup_{n\geq 1}E\mathcal{W}_{1}(\mu_{n}^{N},\mu_{n}) \to 0$ as $N \to \infty$.
\end{enumerate}
\end{Theorem}
\begin{Corollary}\label{cor1}
Suppose Assumptions \ref{assumC} and \ref{assumD} hold and suppose $\delta \in (0, a_0).$  Then
\begin{eqnarray}
\limsup_{N \to\infty}\limsup_{n\to \infty} E \mathcal{W}_{1}(\mu_{n}^{N},\mu_{\infty}) = \limsup_{n\to \infty} \limsup_{N \to\infty}E\mathcal{W}_{1}(\mu_{n}^{N},\mu_{\infty}) =0.\label{3.1}
\end{eqnarray}
\end{Corollary}
The interchangeability of the limits given in Corollary \ref{cor1} allows one to characterize the large $N$ limit of the steady state behavior of the particle system. We need the following assumption on the Markov chain $\{X_{n}(N)\}_{n\geq 0}$ where recall $X_{n}(N):= \{X_{n}^{i},i=1,...,N\}.$ 
The assumption is essentially a communicability condition on the underlying Markovian dynamics.
\begin{Assumption}\label{assumE}
For every $N\geq 1$, the Markov chain  $\{X_{n}(N)\}_{n\geq 0}$ has at most one invariant measure. 
\end{Assumption}
\begin{Theorem}\label{th4}
Suppose  Assumptions \ref{assumA}, \ref{assumB}, \ref{assumC} and \ref{assumE} hold and suppose $\delta \in (0,a_0)$. Then for every $N\in \mathbb{N}$, the Markov chain $\{X_{n}(N)\}_{n\geq 0}$ has a unique invariant measure $\Pi_{\infty}^{N}$. Suppose in addition 
Assumption \ref{assumD} holds. Then $\Pi_{\infty}^{N}$ is $\mu_{\infty}$- chaotic, where $\mu_{\infty}$ is as in Theorem \ref{th2}.
\end{Theorem}
Theorem \ref{th3} gives conditions under which $\mathcal{W}_{1}(\mu_{n}^{N},\mu_{n})$ converges to $0$ as $N \to \infty$, in $L^1$, uniformly in $n$. The next three theorems show that under additional conditions, one can provide concentration bounds uniformly in $n$ which give estimates on the rate of convergence.   Recall the measure $\mu_0$ introduced at the beginning of Section \ref{sec:mod-desc}.

With $\alpha, \sigma_{1}(\alpha)$ defined in Assumption \ref{assumD} and $\omega$ as in Assumption \ref{assumC}, let 
\begin{equation}\label{eq:adefnAB}
a \equiv a(\alpha) := \frac{4^{-\alpha} - e^{-(1+\alpha)\omega}}{2\sigma_{1}(\alpha)}.\end{equation}
\begin{Theorem}\label{th5} 
Suppose Assumptions \ref{assumC} and \ref{assumD} holds.  Fix $\gamma_0 \in (0, a_0)$ and suppose that
 $\delta \in(0,\min\{a^{\frac{1}{1+\alpha}},(a_0-\gamma_0)\})$.  Let $\vartheta = \frac{1-2\sigma \gamma_0}{e^{-\omega}+ 2\delta\sigma}$.
Then there exists $N_0 \in \mathbb{N}_0$ and $C_1 \in (0, \infty)$ such that for all $\eps > 0,$ and for all $n\ge 0,$
$$ P(\mathcal{W}_{1}(\mu_{n}^{N},\mu_{n})>\varepsilon) \leq P(\mathcal{W}_{1}(\mu_{0}^{N},\mu_{0})>2\sigma\gamma_0\vartheta^n\varepsilon) + C_1 \varepsilon^{-(1+\alpha)} N^{-\frac{\alpha}{d+2}}, $$
for all  $N >N_{0}
\left(\max\left\{1,\log^+ \varepsilon\right\}\right)^{\frac{d+2}{d}}$.
\end{Theorem}
\begin{remark}\label{Cor2}
	\begin{description}
		\item{(i)}  Since $\delta < a_0 - \gamma_0$, we have that $\vartheta > 1$ and so the above theorem gives the following uniform concentration estimate: 
		$$ \sup_{n \ge 1} P(\mathcal{W}_{1}(\mu_{n}^{N},\mu_{n})>\varepsilon) \leq P(\mathcal{W}_{1}(\mu_{0}^{N},\mu_{0})>2\sigma\gamma_0\varepsilon) + C_1 \varepsilon^{-(1+\alpha)} N^{-\frac{\alpha}{d+2}}, $$
		for all  $N >N_{0}
		\left(\max\left\{1,\log^+ \varepsilon\right\}\right)^{\frac{d+2}{d}}$.
\item{(ii)} Under additional conditions on $\{X_0^{i,N}\}$ one can give concentration bounds for the first term on the right side of the above inequality.  For example, when
$\{X_0^{i,N}\}_{i=1}^N$ are i.i.d. such concentration bounds can be found in Theorem 2.7 of \cite{bolley2007quantitative}.  Also, although not pursued here, the bound obtained in Theorem \ref{th5}
can be used to give an alternative proof of Theorem \ref{th3}(2).
\end{description}
\end{remark}
Next we obtain exponential concentration bounds.  The bounds depend in particular on our assumptions on the initial condition.  Our first result (Theorem \ref{th6}) treats the case
where the initial random vector has a general exchangeable distribution while the second result  (Theorem \ref{thm6}) considers a more restrictive setting where the initial random vector is i.i.d.  In the second case the probabilities will decay exponentially in $N$ whereas in the first case the exponent will be some dimensional dependent power of $N$. 

We start with our main assumption for Theorem \ref{th6}.
\begin{Assumption}\label{assumG}
	\begin{description}
		\item{(i)} For some $M \in (1, \infty)$, $D(x) \le M$ for $\theta$ a.e. $x \in \R^m$. 
	\item{(ii)} There exists $\alpha \in (0, \infty)$ such that 
 $\int e^{\alpha |x|} \mu_0(dx) < \infty$ and  $\int e^{\alpha D_1(z)} \theta(dz) < \infty$.
		\end{description}
\end{Assumption}
\begin{Theorem}\label{th6} Suppose that Assumptions  \ref{assumC}  and \ref{assumG} hold.  Fix $\gamma_0 \in (0, a_0)$ and suppose that
	$\delta \in [0,\min\{a_0-\gamma_0,\frac{1-e^{\omega}}{2M}\})$. Then there exists $N_0 \in \mathbb{N}$ and $C_1 \in (0, \infty)$ such that for all
	$\e > 0$ 
	$$P[\mathcal{W}_{1}(\mu_{n}^{N},\mu_{n})>\varepsilon]\leq P[\mathcal{W}_{1}(\mu_{0}^{N},\mu_{0})>2\sigma\gamma_0 \vartheta^n\varepsilon]
	+ e^{-C_1 \e N^{1/d+2}},$$
	for all  $n \ge 0$, $N \ge N_0 \max\{ (\frac{1}{\e} \log^+ \frac{1}{\e})^{d+2}, \e^{(d+2)/(d-1)}\}$, if $d > 1$; and
	$$P[\mathcal{W}_{1}(\mu_{n}^{N},\mu_{n})>\varepsilon]\leq P[\mathcal{W}_{1}(\mu_{0}^{N},\mu_{0})>2\sigma\gamma_0 \vartheta^n\varepsilon]
	+ e^{-C_1 (\e \wedge 1) N^{1/d+2}},$$
	for all $n\ge 0$, $N \ge N_0 \max\{ (\frac{1}{\e} \log^+ \frac{1}{\e})^{d+2}, 1\}$, if $d = 1$.  Here $\vartheta \in (1, \infty)$ is as in Theorem \ref{th5}.
	\end{Theorem}
%
Finally we consider the case where the initial distribution of the $N$ particles is i.i.d.  The proof relies on various estimates from \cite{boissard2011simple, bolley2005weighted}.
\begin{Theorem}\label{thm6}
Suppose that $\{X_{0}^{i,N}\}_{i=1,..,N}$ are i.i.d. with common distribution  $\mu_{0}$ for each $N$. Suppose that Assumptions  \ref{assumC} and \ref{assumG} hold.
Fix $\gamma \in (0,1-e^{-\omega})$.
 Suppose that  
$\delta \in \left[0, \frac{1 -e^{-\omega}-\gamma}{2M}\right).$
Then 
there exist $N_0, a_1, a_2 \in (0, \infty)$ and a nonincreasing function $\vs_1: (0, \infty) \to (0,\infty)$ such that $\vs_1(t) \downarrow 0$ as $t \uparrow \infty$ and
for all $\e > 0$ and $N \ge N_0 \vs_1(\e)$
\[\sup_{n\ge 0}P[\mathcal{W}_{1}(\mu_{n}^{N},\mu_{n})>\varepsilon] \le 
  a_1e^{-N a_2 (\e^2 \wedge \e)} 
\]
\end{Theorem}
\begin{remark}
	\label{rem:vsconstant}
	(i) One can describe the function $\vs_1$ in the above theorem explicitly.  
	Define for $\gamma \in (0,1)$, $m_{\gamma}: (0, \infty) \to (0, \infty)$ as $m_{\gamma}(t) = \frac{\gamma t}{\delta M}$, where $M$ is as in Assumption \ref{assumG}.
	Then
	$$\vs_1(t) = \max \left\{1, \frac{\log \mathcal{C}^0_{m_{\gamma}(t)}}{m_{\gamma}^2(t)}, \frac{\log \mathcal{C}^0_{\gamma t}}{\gamma^2t^2}, \frac{1}{t^2}, \frac{1}{t}\right\},$$
	where $\mathcal{C}^0_t$ is defined by the right side of \eqref{eq:eq549} with $\zeta$ replaced by $\zeta_0$ where $\zeta_0$ is as in Corollary \ref{cor:cor621}.\\ \ \\
	
	\noindent (ii)  If Assumption \ref{assumG} is strengthened to $\int e^{\alpha D_1(z)^2} \theta(dz) < \infty$ for some $\alpha > 0$ then one can 
	strengthen the conclusion of Theorem \ref{thm6} as follows:
	For $\delta$ sufficiently small
	there exist $N_0, a_1, a_2 \in (0, \infty)$ and a nonincreasing function $\vs_2: (0, \infty) \to (0,\infty)$ such that $\vs_2(t) \downarrow 0$ as $t \uparrow \infty$ and
	for all $\e > 0$ and $N \ge N_0 \vs_2(\e)$
	\[\sup_{n\ge 0}P[\mathcal{W}_{1}(\mu_{n}^{N},\mu_{n})>\varepsilon] \le 
	    a_1e^{-N a_2 \e^2 }. 
	\]
	\end{remark}
\section{Proofs}
The following elementary lemma gives a basic moment bound that will be used in our analysis.
\subsection{Proof of Theorem \ref{th1}}
\begin{lemma}\label{lem1}
 Suppose Assumptions \ref{assumA} and \ref{assumB} hold. Then, for every $n\geq 1$, 
$$M_{n} := \sup_{N\geq 1} \max_{1\leq i \leq N} E |X_{n}^{i}| < \infty .$$
In addition, if Assumption \ref{assumC} holds and $\delta \in (0,a_0)$ then $\sup_{n\geq 1} M_{n} <\infty.$
\end{lemma}
\textbf{Proof:} 
We will only prove the second statement in the lemma. Proof of the first statement is similar.
Note that, for $n\geq 1$ and $i=1,..,N$
$$X_{n}^{i} =AX_{n-1}^{i}+\delta f(X_{n-1}^{i},\mu_{n-1}^{N},\epsilon_{n}^{i}).$$
Thus 
\begin{eqnarray}\label{eq1.1.1}
|X_{n}^{i}| \leq ||A||\ |X_{n-1}^{i}|+ \delta D(\epsilon_{n}^{i}) [|X_{n-1}^{i}| + ||\mu_{n-1}^{N}||_{1}] + \delta D_{1}(\epsilon_{n}^{i}).
\end{eqnarray}
From exchangeability of  $\{X_{n-1}^{k},k=1,...,N\}$ it follows that 
$$E||\mu_{n-1}^{N}||_{1} = E[\int |x| d\mu_{n-1}^{N}] = E\frac{1}{N}\sum_{k=1}^{N}|X_{n-1}^{k}|= E|X_{n-1}^{1}|.$$
Taking expectation  in (\ref{eq1.1.1}) and using independence between $\epsilon_{n}^{i}$ and $\{X_{n-1}^{j}\}_{j= 1}^{N}$, we have
\begin{eqnarray}
E|X_{n}^{i}| \leq (||A||+ 2\delta \sigma) E|X_{n-1}^{i}| + \delta c_{0}.\label{eq1.1.2}
\end{eqnarray}
The assumption on $\delta$ implies that $\gamma:= ||A||+2\delta \sigma \in (0,1) .$
A recursive application of (\ref{eq1.1.2}) now shows that
$$E|X_{n}^{i}| \leq \gamma^{n} E|X_{0}^{i}| + \frac{\delta c_{0}}{1-\gamma}.$$
The result follows.
\qed

Recall the map $\Psi$ defined in \eqref{eq:544}.
\begin{lemma}\label{cor2}
Under Assumptions \ref{assumA} and \ref{assumB},  for every $\epsilon >0$ and $n\geq 1$, there exists a compact set $K_{\epsilon,n} \in \mathcal{B}(\mathbb{R}^{d})$ such that 
$$\sup_{N\geq 1} E\left\{ \int_{K_{\epsilon,n}^{c}} |x| \left(\mu_{n}^{N}(dx)+ \Psi(\mu_{n-1}^{N})(dx)\right)\right\}<\epsilon .$$
\end{lemma}
\textbf{Proof:} Note that for any non-negative $\phi: \mathbb{R}^{d} \to \mathbb{R}$,  
\begin{eqnarray}
E \int \phi(x) \mu_{n}^{N}(dx) = \frac{1}{N} \sum_{k=1}^{N} E \phi(X_{n}^{k}) = E\phi(X_{n}^{1}),
\label{sum}  
\end{eqnarray}
and
\begin{eqnarray} 
E\int \phi(x) \Psi(\mu_{n}^{N})(dx) &=& \frac{1}{N} \sum_{i=1}^{N} E(E(\langle \phi , \delta_{X_{n}^{i}}P^{\mu_{n}^{N}}\rangle \mid \mathcal{F}_{n}))\nonumber\\
&=& \frac{1}{N}\sum_{i=1}^{N} E\phi\left(AX_{n}^{i} + \delta f(X_{n}^{i}, \mu_{n}^{N}, \epsilon_{n+1}^{i})\right)\nonumber\\
&=&\frac{1}{N}\sum_{i=1}^{N} E\phi(X_{n+1}^{i})=  E\phi(X_{n+1}^{1}). \label{eq1.2.1}
\end{eqnarray}
To get the desired result from above equalities it suffices to show that 
\begin{eqnarray}
\text{the family $\{X_{n}^{i,N},i=1,...,N; N\geq 1\}$ is uniformly integrable for every $n\geq 0$.}\label{e1} 
\end{eqnarray}
We will prove (\ref{e1}) by induction on $n$. 
Once more we suppress $N$ from the super-script. 
Clearly by our assumptions $\{X_{0}^{i},i=1,...,N; N\geq 1\}$ is uniformly integrable. Now suppose that the Statement (\ref{e1}) holds for some $n$. Note that
\begin{eqnarray}
|X_{n+1}^{i}| &\leq& ||A||\ |X_{n}^{i}|+ \delta D(\epsilon_{n+1}^{i}) [|X_{n}^{i}| + ||\mu_{n}^{N}||_{1}] + \delta D_{1}(\epsilon_{n+1}^{i})\nonumber\\
&=& ||A||\ |X_{n}^{i}|+ \delta D(\epsilon_{n+1}^{i}) [ |X_{n}^{i}| + \frac{1}{N}\sum_{i=1}^{N}|X_{n}^{i}|] + \delta D_{1}(\epsilon_{n+1}^{i}).\label{eq1.2.2}
\end{eqnarray}
From exchangeability it follows $\frac{1}{N}\sum_{i=1}^{N}|X_{n}^{i}| = E[|X_{n}^{i}| \mid \mathcal{G}_{n}^{N}]$, where $\mathcal{G}_{n}^{N} = \sigma\{\frac{1}{N}\sum_{i=1}^{N}\delta_{X_{n}^{i}}\}$. Combining this with the induction hypothesis that $\{X_{n}^{i},i=1,...,N; N\geq 1\}$ is uniformly integrable, we see that $\{\frac{1}{N}\sum_{i=1}^{N}|X_{n}^{i}|,N\geq 1\}$ is uniformly integrable. Here we have used the fact that if $\{Z_{\alpha}, \alpha \in \Gamma_1\}$ is a uniformly integrable family
and $\{\clh_{\beta}, \beta \in \Gamma_2\}$ is a collection of $\sigma$-fields where $\Gamma_1, \Gamma_2$ are arbitrary index sets, then
$\{E (Z_{\alpha} \mid \clh_{\beta}), (\alpha, \beta) \in \Gamma_1\times \Gamma_2\}$ is a uniformly integrable family.
Also from Assumptions \ref{assumA} and \ref{assumB} the families $\{D(\epsilon_{n+1}^{i});i\geq 1\}$, $\{D_{1}(\epsilon_{n+1}^{i});i\geq 1\}$ are uniformly integrable. These observations along with independence between $\{\epsilon_{n+1}^{i},i=1,..,N\}$ and $\{X_{n}^{i}:i=1,...,N;N\geq 1\}$ yield that the family $\{|X_{n}^{i}|:i=1,...,N;N\geq 1\}$ is uniformly integrable. The result follows.
\qed

We now proceed to the proof of Theorem \ref{th1}. 
We will argue via induction on $n\geq 0$. By  assumption (\ref{pointwise}) holds for $n=0$. Assume now that it holds for some $n>0$ . Note that,
\begin{eqnarray}
\mathcal{W}_{1}(\mu_{n+1}^{N},\mu_{n+1})  &\leq& \mathcal{W}_{1}(\mu_{n+1}^{N},\mu_{n}^{N} P^{ \mu_{n}^{N}})  +\mathcal{W}_{1}(\mu_{n}^{N} P^{ \mu_{n}^{N}}, \mu_{n}^{N} P^{ \mu_{n}}) +\mathcal{W}_{1}(\mu_{n}^{N} P^{ \mu_{n}},\mu_{n+1}).\quad  \label{eq1.3.1}
\end{eqnarray}
Consider the last term in (\ref{eq1.3.1}).  Using Assumption \ref{assumA} we see that if $\phi$ is Lipschitz then $ P^{\mu_{n}}\phi$ is Lipschitz and $||P^{\mu_{n}}\phi||_{1}\leq (||A||+\delta \sigma )||\phi||_{1}.$ Thus, almost surely
\begin{eqnarray}
\sup_{\phi \in \mbox{Lip}_{1}}|\langle \phi,\mu_{n}^{N} P^{ \mu_{n}}-\mu_{n+1}\rangle|&=&\sup_{\phi \in \mbox{Lip}_{1}}|\langle P^{\mu_{n}}\phi,\mu_{n}^{N} -\mu_{n}\rangle| \nonumber\\ &\leq& (||A||+\delta \sigma )    \sup_{g \in \mbox{Lip}_{1}}|\langle g, \mu_{n}^{N} -\mu_{n}\rangle|\quad .\nonumber 
\end{eqnarray}
Taking expectations we obtain,
\begin{eqnarray}
E \mathcal{W}_{1}(\mu_{n}^{N}P^{ \mu_{n}},\mu_{n+1}) \leq (||A||+\delta \sigma) E \mathcal{W}_{1}(\mu_{n}^{N},\mu_{n}).\label{eq1.3.2}
\end{eqnarray}
Consider now the second term in (\ref{eq1.3.1}). Using Assumption \ref{assumA} again, we have,
\begin{eqnarray}
\sup_{\phi\in \mbox{Lip}_{1}}|\langle \phi,\mu_{n}^{N} P^{ \mu_{n}^{N}} - \mu_{n}^{N} P^{ \mu_{n}}\rangle |&\leq& \frac{1}{N}\sum_{i=1}^{N}\int \left[|\phi(AX_{n}^{i} + \delta f(X_{n}^{i},\mu_{n}^{N},\xi)) \right. \nonumber\\&-& \left . \phi(AX_{n}^{i} +\delta f(X_{n}^{i},\mu_{n},\xi))|\right]\theta(d\xi)\nonumber\\
&\leq&  \delta \sigma\mathcal{W}_{1}(\mu_{n}^{N},\mu_{n}).\nonumber
\end{eqnarray}
Taking expectations we get
\begin{eqnarray}
E\mathcal{W}_{1}(\mu_{n}^{N}P^{\mu_{n}^{N}},\mu_{n}^{N}P^{\mu_{n}})= E\sup_{\phi\in \mbox{Lip}_{1}(\mathbb{R}^{d})}|\langle \phi,\mu_{n}^{N} P^{ \mu_{n}^{N}} - \mu_{n}^{N} P^{ \mu_{n}}\rangle| &\leq&\delta \sigma E \mathcal{W}_{1}(\mu_{n}^{N},\mu_{n}).\label{eq:ab428}
\end{eqnarray}
Now we consider the first term of the right hand side of (\ref{eq1.3.1}). We will use Lemma \ref{cor2}. Fix $\epsilon>0$ and let $K_{\epsilon}$ be a compact set in $\mathbb{R}^{d}$ such that
$$\sup_{N\geq 1} E\left\{ \int_{K_{\epsilon}^{c}} |x|(\mu_{n+1}^{N}(dx)+ \Psi(\mu_{n}^{N})(dx))\right\}<\epsilon .$$
Let $\Lip_{1}^{0}(\mathbb{R}^{d}):= \{f\in \Lip_{1}(\mathbb{R}^{d}): f(0)=0\}$.
Then,
\begin{eqnarray}
E\sup_{\phi \in \mbox{Lip}_{1}(\mathbb{R}^{d})}|\langle\phi, \mu_{n+1}^{N} - \mu_{n}^{N} P^{\mu_{n}^{N}}\rangle| &=& E\sup_{\phi \in \Lip_{1}^{0}(\mathbb{R}^{d})}|\langle\phi, \mu_{n+1}^{N} - \mu_{n}^{N} P^{\mu_{n}^{N}}\rangle|\nonumber\\
&\leq& E\sup_{\phi \in \Lip_{1}^{0}(\mathbb{R}^{d})}|\langle\phi. 1_{K_{\epsilon}}, \mu_{n+1}^{N} - \mu_{n}^{N} P^{\mu_{n}^{N}}\rangle| +\epsilon.\label{1.3.4} 
\end{eqnarray}
We will now apply Lemma \ref{app1} in the Appendix. Note that for any $\phi\in \Lip_{1}^{0}(\mathbb{R}^{d})$, $\sup_{x \in K_{\epsilon}} |\phi(x)|\leq diam(K_{\epsilon}) := m_{\epsilon}.$

Thus with notation as in Lemma \ref{app1} 
\begin{eqnarray}
 \sup_{\phi \in \Lip_{1}^{0}(\mathbb{R}^{d})}|\langle\phi. 1_{K_{\epsilon}}, \mu_{n+1}^{N} - \mu_{n}^{N} P^{\mu_{n}^{N}}\rangle| 
&\leq& \max_{\phi \in \mathcal{F}^{\epsilon}_{m_{\epsilon,1}}(K_{\epsilon})}|\langle\phi , \mu_{n+1}^{N} - \mu_{n}^{N} P^{\mu_{n}^{N}}\rangle| +2\epsilon.\label{1.3.5}
\end{eqnarray}
where we have denoted the restrictions of $\mu_{n+1}^{N}$ and $\mu_{n}^{N}P^{\mu_{n}}$ to $K_{\epsilon}$ by the same symbols.
Using the above inequality in (\ref{1.3.4}), we obtain 
\begin{eqnarray}
E\mathcal{W}_{1}( \mu_{n+1}^{N} , \mu_{n}^{N} P^{\mu_{n}^{N}})  &\leq& \sum_{\phi \in \mathcal{F}^{\epsilon}_{m_{\epsilon,1}}(K_{\epsilon})} E|\langle \phi, \mu_{n+1}^{N} - \mu_{n}^{N} P^{\mu_{n}^{N}}  \rangle| +3\epsilon.\nonumber
\end{eqnarray}
Using Lemma \ref{app2} we see that the second term on the left hand side can be bounded by $\frac{2m_{\epsilon}| \mathcal{F}^{\epsilon}_{m_{\epsilon,1}}(K_{\epsilon})|}{\sqrt{N}}$. Combining this estimate with \eqref{eq1.3.1}, (\ref{eq1.3.2}) and (\ref{eq:ab428}) we now have
\begin{eqnarray}
E\mathcal{W}_{1}(\mu_{n+1}^{N},\mu_{n+1}) &\leq& (||A||+ 2\delta \sigma) E\mathcal{W}_{1}(\mu_{n}^{N},\mu_{n})  + \frac{2m_{\epsilon}| \mathcal{F}^{\epsilon}_{m_{\epsilon,1}}(K_{\epsilon})|}{\sqrt{N}} +3\epsilon. \label{1.3.6}
\end{eqnarray}
Sending $N\to \infty$ in (\ref{1.3.6}) and using induction hypothesis, we have 
$$\limsup_{N\to \infty}E\mathcal{W}_{1}(\mu_{n+1}^{N},\mu_{n+1}) \leq 3\epsilon.$$
Since $\epsilon >0$ is arbitrary, the result follows.
\qed

\label{sec-int}

\subsection{Proof of Theorem \ref{th2}}
We begin with the following estimate.
\begin{lemma}\label{lem2.1}
Under Assumptions \ref{assumA},\ref{assumB} and \ref{assumC}
$$\mathcal{W}_{1}(\Psi^{n}(\mu_{0}),\Psi^{n}(\gamma_{0})) \leq (e^{-\omega} +2\delta \sigma)^{n}\ \mathcal{W}_{1}(\mu_{0},\gamma_{0})$$
for any choice of $\mu_{0},\gamma_{0} \in \mathcal{P}_{1}(\mathbb{R}^{d})$.
\end{lemma}
\textbf{Proof:} Given $\mu_{0},\gamma_{0} \in \mathcal{P}_{1}(\mathbb{R}^{d}),$ let $\mathcal{C}(\mu_{0},\gamma_{0})=\{\mu\in \mathcal{P}_{1}(\mathbb{R}^{d}\times\mathbb{R}^{d})\mid  \mu_{0}(\cdot)=\mu(\cdot \times \mathbb{R}^{d}),\gamma_{0}(\cdot)=\mu(\mathbb{R}^{d}\times \cdot) \}.$
Fix $\mu \in \mathcal{C}(\mu_{0},\gamma_{0})$ and let $(X_{0},Y_{0})$ be $\mathbb{R}^{d}\times \mathbb{R}^{d}$ valued random variables with distribution $\mu$. Also, let $\{\epsilon_{n}\}_{n\geq 1}$ be an iid sequence of random variables with common law $\theta$ independent of $(X_{0},Y_{0}).$ Define for $n\geq 0$, 
\begin{eqnarray}
X_{n+1} &=& AX_{n} +\delta f(X_{n},\mu_{n},\epsilon_{n+1}),\nonumber\\
Y_{n+1}&=& AY_{n} + \delta f(Y_{n},\gamma_{n},\epsilon_{n+1}) \nonumber
\end{eqnarray}
where $\mu_{n} = \mathcal{L}(X_{n})$ and $\gamma_{n} = \mathcal{L}(Y_{n})$. Then clearly $\mu_{n}=\Psi^{n}(\mu_{0}),\gamma_{n}=\Psi^{n}(\gamma_{0})$. For $n\geq 0,$ denote $\beta_{n}=\mathcal{W}_{1}(\mu_{n},\gamma_{n}),\alpha_{n}= E|X_{n}-Y_{n}|.$ Then
\begin{eqnarray}
\beta_{n+1}&=& \sup_{\phi \in \mbox{Lip}_1}\{|\int \phi d\mu_{n+1}- \int\phi d\gamma_{n+1}|\}\nonumber\\
&=& \sup_{\phi\in \mbox{Lip}_{1}}\{|E \phi (X_{n+1})- E\phi(Y_{n+1})| \}\nonumber\\
&\leq& E|X_{n+1}-Y_{n+1}| = \alpha_{n+1}. \hspace{2cm}\label{eq2.1}
\end{eqnarray}
Also,
\begin{eqnarray}
\alpha_{n+1}&\leq& ||A||E|X_{n}-Y_{n}|+\delta E|f(X_{n},\mu_{n},\epsilon_{n+1})-f(Y_{n},\gamma_{n},\epsilon_{n+1})|   \nonumber\\
&\leq& e^{-\omega}E|X_{n}-Y_{n}| +\delta \sigma (E|X_{n}-Y_{n}| + \mathcal{W}_{1}(\mu_{n},\gamma_{n}))\label{eq2.2}\hspace{1cm}\nonumber\\
&=& (e^{-\omega}+\delta \sigma)\alpha_{n} +\delta \sigma\beta_{n}, \nonumber\\
&\leq & (e^{-\omega} + 2\delta \sigma)\alpha_{n}\label{eq2.3}
\end{eqnarray}
where the second inequality in the display follows from  Assumptions \ref{assumA} and \ref{assumC}. Combining \eqref{eq2.1} and \eqref{eq2.3} we have
\begin{eqnarray} 
\beta_{n+1} &\leq& (e^{-\omega}+2\delta \sigma)^{n+1} E|X_{0} - Y_{0}|\nonumber\\
&=&  (e^{-\omega}+2\delta \sigma)^{n+1} \int|x-y|\mu(dx dy).\nonumber
\end{eqnarray}
 We now have, on taking infimum on the right hand side of the above display over all
 $\mu \in \mathcal{C}(\mu_{0},\gamma_{0})$, that
 $\beta_{n+1} \leq (e^{-\omega} + 2\delta \sigma)^{n+1}\beta_{0}$. The result follows.
\qed

We now complete the proof of Theorem \ref{th2}.  Observe that under our assumption on $\delta$, $\chi := e^{-\omega}+2\delta \sigma \in (0,1)$.
The first part of the theorem now follows from 
 Lemma \ref{lem2.1} and Banach's fixed point theorem.
Furthermore
\begin{eqnarray}
\mathcal{W}_{1}(\Psi^{n}(\mu),\mu_{\infty}) = \mathcal{W}_{1}(\Psi^{n}(\mu),\Psi^{n}(\mu_{\infty}))
\leq  \chi^{n}\ \mathcal{W}_{1}(\mu,\mu_{\infty}).
\end{eqnarray}
Second part of the theorem is now immediate. 
\qed

\subsection{Proof of Theorem \ref{th3}}
We start with the following moment bound.
\begin{lemma}\label{lem2}
Suppose that Assumptions \ref{assumC} and \ref{assumD} hold and suppose that 
$\delta \in (0,a(\alpha)^{\frac{1}{1+\alpha}})$ where $a(\cdot)$ is as in \eqref{eq:adefnAB}. Then $\sup_{N \ge 1}\sup_{n \geq 1} E|X_{n}^{1,N}|^{1+\alpha} < \infty.$
\end{lemma}
\textbf{Proof:} By Holder's inequality,  for any four nonnegative real numbers $a_{1},a_{2},a_{3},a_{4},$ 
\begin{eqnarray}
(a_{1}+a_{2}+a_{3}+a_{4})^{1+\alpha} \leq 4^{\alpha}(a_{1}^{1+\alpha}+a_{2}^{1+\alpha}+a_{3}^{1+\alpha}+a_{4}^{1+\alpha}).\label{4.3.1} 
\end{eqnarray}
Using Assumption \ref{assumC} and (\ref{fbound})
\begin{eqnarray}
|X_{n+1}^{i}| &\leq& e^{-\omega}|X_{n}^{i}| + \delta |D(\epsilon_{n+1}^{N})| \left(|X_{n}^{i}| + |\mu_{n}^{N}|\right) + D_{1}(\epsilon_{n+1}^{N}).\nonumber 
\end{eqnarray}
Taking expectations on both sides and applying (\ref{4.3.1}), we have, from Assumption \ref{assumD}
\begin{eqnarray}
E|X_{n+1}^{i}|^{1+\alpha}  &\leq& 4^{\alpha}e^{-\omega(1+\alpha)} E|X_{n}^{i}|^{1+\alpha} +4^{\alpha} \delta^{(1+\alpha)}\sigma_{1}[E|X_{n}^{i}|^{1+\alpha} + E|\mu_{n}^{N}|^{1+\alpha}]\nonumber +4^{\alpha}c_{1}\nonumber\\
&\leq& 4^{\alpha}e^{-\omega(1+\alpha)} E|X_{n}^{i}|^{1+\alpha} +4^{\alpha} \delta^{(1+\alpha)}\sigma_{1}[2E|X_{n}^{i}|^{1+\alpha} ]+4^{\alpha}c_{1}\quad \label{4.3.1.1}
\end{eqnarray}
where the last line in the display follows from Jensen's inequality: $ E|\mu_{n}^{N}|^{1+\alpha} =  E|\int |x| \mu_{n}^{N}(dx)|^{1+\alpha} \leq E|\int |x| ^{1+\alpha}\mu_{n}^{N}(dx)| =  E|X_{n}^{i}|^{1+\alpha} $ and $c_{1} $ is as in Assumption \ref{assumD}. 

Note that under our  condition on $\delta$ 
$$\kappa_1 \equiv 4^{\alpha}[e^{-\omega(1+\alpha)} + 2\delta^{(1+\alpha)}\sigma_{1}] <1.$$
Thus
\begin{equation}\sup_{n\ge 1} E|X_n^i|^{1+\alpha} \le \kappa_1 E|X_0^i|^{1+\alpha} + \frac{\kappa_2}{1-\kappa_1}, \label{eq:eq258}
\end{equation}
where $\kappa^2 = 4^{\alpha} c_1$. The result follows. 
\qed

We now complete the proof of Theorem \ref{th3}. Note that
$$\mu_{n}^{N} -\mu_{n} = \sum_{k=1}^{n} [\Psi^{n-k}(\mu_{k}^{N}) - \Psi^{n-k}.\Psi(\mu_{k-1}^{N})] + [\Psi^{n}(\mu_{0}^{N}) - \Psi^{n}(\mu_{0})]. $$
It then follows using Lemma \ref{lem2.1} that with $\chi = (e^{-\omega} + 2\delta \sigma),$ almost surely
\begin{eqnarray}
\mathcal{W}_{1}(\mu_{n}^{N} , \mu_{n})&\leq& \sum_{k=1}^{n}\mathcal{W}_{1}(\Psi^{n-k}(\mu_{k}^{N}) ,  \Psi^{n-k}(\Psi(\mu_{k-1}^{N})))+\mathcal{W}_{1}(\Psi^{n}(\mu_{0}^{N}), \Psi^{n}(\mu_{0}))\nonumber\\
&\leq& \sum_{k=1}^{n}\chi^{n-k}\ \mathcal{W}_{1}(\mu_{k}^{N},\Psi(\mu_{k-1}^{N})) +\chi^{n}\mathcal{W}_{1}(\mu_{0}^{N},\mu_{0}).  \label{4.3.2}
\end{eqnarray}
Taking expectations,
\begin{eqnarray}
E\mathcal{W}_{1}(\mu_{n}^{N},\mu_{n}) &\leq& \sum_{k=1}^{n} \chi ^{n-k} \ E\mathcal{W}_{1}(\mu_{k}^{N},\Psi(\mu_{k-1}^{N})) + \chi^{n}\ E\mathcal{W}_{1}(\mu_{0}^{N},\mu_{0}).\nonumber 
\end{eqnarray}
Since $a(\alpha_0)^{1/(1+\alpha_0)} \to a_0$ as $\alpha_0 \to 0$ and $\delta \in (0, a_0)$, we can find $\alpha_0  \in (0, \alpha)$ such that 
$\delta \in (0, a(\alpha_0)^{1/(1+\alpha_0)})$.  From Lemma \ref{lem2} we then have that $\sup_{N \ge 1}\sup_{n \geq 1} E|X_{n}^{1}|^{1+\alpha_0} < \infty$
and consequently  the family $\{ X_{n}^{i}, i=1,...,N, N\geq 1, n\geq 1\}$ is uniformly integrable. 
Similar to the proof of Corollary \ref{cor2} (cf. the argument below \eqref{eq1.2.2}) using
 (\ref{sum}) and (\ref{eq1.2.1}) it follows that, for some compact $K_{\epsilon} \subseteq \mathbb{R}^{d} $
\begin{eqnarray}
\sup_{N\geq 1} \sup_{n\geq 1} E \int_{K_{\epsilon}^{c}} |x| \ [\mu_{n}^{N}(dx) + \Psi(\mu_{n-1}^{N})(dx)] < \epsilon. \label{4.3.3}
\end{eqnarray}
Now for every $k \ge 1$ 
\begin{eqnarray}
E\mathcal{W}_{1}(\mu_{k}^{N},\Psi(\mu_{k-1}^{N})) &=&  E\sup_{f\in \Lip_{1}^{0}(\mathbb{R}^{d})}|\langle f ,\mu_{k}^{N} - \Psi(\mu_{k-1}^{N})\rangle | \nonumber\\&\leq& E \sup_{f\in \Lip_{1}^{0}(\mathbb{R}^{d})}|\langle f.1_{K_{\epsilon}} ,\mu_{k}^{N} - \Psi(\mu_{k-1}^{N})\rangle | + \epsilon, \nonumber
\end{eqnarray}
when $Lip_{1}^{0}(\mathbb{R}^{d})$ is as introduced above (\ref{1.3.4}). 
Applying Lemmas \ref{app1} and \ref{app2} as in the proof of Theorem \ref{th1} we now see that
$$E\mathcal{W}_{1}(\mu_{k}^{N},\Psi(\mu_{k-1}^{N})) \leq |\mathcal{F}^{\epsilon}_{m_{\epsilon},1}(K_{\epsilon})|\frac{2 m_{\epsilon}}{\sqrt{N}} + 3\epsilon,$$
where $m_{\epsilon} = diam(K_{\epsilon})$.
Thus 
\begin{eqnarray}
E\mathcal{W}_{1}(\mu_{n}^{N},\mu_{n}) &\leq& \sum_{k=1}^{n} \chi^{n-k}\ \{|\mathcal{F}^{\epsilon}_{m_{\epsilon},1}(K_{\epsilon})|\frac{2m_{\epsilon}}{\sqrt{N}} + 3\epsilon\} + \chi^{n} E\mathcal{W}_{1} (\mu_{0}^{N},\mu_{0})\nonumber \\
&\leq& \{|\mathcal{F}^{\epsilon}_{m_{\epsilon},1}(K_{\epsilon})|\frac{2m_{\epsilon}}{\sqrt{N}} + 2\epsilon\} \frac{1}{1-\chi} + \chi^{n}\ E\mathcal{W}_{1}(\mu_{0}^{N},\mu_{0}). \label{eq 3.1}
\end{eqnarray}
Given $\varepsilon >0,$ choose $\epsilon$ sufficiently small and $N_{0}$ sufficiently large such that $\forall N\geq N_{0}$ 
$$\{|\mathcal{F}^{\epsilon}_{m_{\epsilon},1}(K_{\epsilon})|\frac{2m_{\epsilon}}{\sqrt{N}} + 2\epsilon\} \frac{1}{1-\chi} \leq \frac{\varepsilon}{2}.$$
Choose $n_{0}$ large enough so that $\forall n\geq n_{0}$, $2\chi^{n} \ ||\mu_{0}||_{1} < \frac{\varepsilon}{2}$. Combining the above estimates we have $\forall N\geq N_{0},$ and $n\geq n_{0}$ 
$E\mathcal{W}_{1}(\mu_{n}^{N},\mu_{n}) \leq \varepsilon.$
This proves the first part of the theorem.

Second part is immediate from the first part and  Theorem \ref{th1}.
\qed

\subsection{Proof of Corollary \ref{cor1}}
\textbf{Proof:} 
Note that
$$
E \mathcal{W}_{1}(\mu_{n}^{N},\mu_{\infty}) \le E \mathcal{W}_{1}(\mu_{n}^{N},\mu_{n}) +  \mathcal{W}_{1}(\mu_{n},\mu_{\infty}).$$
Combining this with \eqref{eq 3.1} we have 
%
%
%
%
\begin{eqnarray}\label{4.4.1}
E \mathcal{W}_{1}(\mu_{n}^{N},\mu_{\infty})  \leq \left(|\mathcal{F}^{\epsilon}_{m_{\epsilon},1}(K_{\epsilon})|\frac{2m_{\epsilon}}{\sqrt{N}} + 2\epsilon\right) \frac{1}{1-\chi} + \chi^{n}\ E\mathcal{W}_{1}(\mu_{0}^{N},\mu_{0}) + \mathcal{W}_{1}(\mu_n,\mu_{\infty}).
\end{eqnarray}
The result now follows on using Theorem \ref{th2}. \qed

%

\subsection{Proof of Theorem \ref{th4}}
For $N\geq 1$ and $n\in \mathbb{N}_{0}$, define $\Pi_{n}^{N} \in \mathcal{P}((\mathbb{R}^{d})^{N})$ as 
\begin{eqnarray}
\langle \phi,\Pi_{n}^{N}  \rangle &=& \frac{1}{n} \sum_{j=1}^{n} E \phi(X_{j}^{1,N},...,X_{j}^{N,N}), \quad \phi\in BM((\mathbb{R}^{d})^{N})
\end{eqnarray}
where $\{X_{j}^{i,N}, j\in \mathbb{N}_{0},i=1,..,N\}$ are as defined in (\ref{nlps}).

From Lemma \ref{lem1} it folows that, for each $N\geq 1,$ the sequence $\{\Pi_{n}^{N},n\geq 1\}$ is relatively compact and using Assumption \ref{assumA} 
it is easy to see that any limit point of $\Pi_n^N$ (as $n\to \infty$) is an invariant measure of the Markov chain $\{X_n(N)\}_{n\ge 0}$. Combining this with
Assumption \ref{assumE} we have that  there exists a unique invariant measure $\Pi_{\infty}^{N} \in \mathcal{P}((\mathbb{R}^{d})^{N})$ for this Markov chain and, as $n \to \infty$,
\begin{eqnarray}
\Pi_{n}^{N} \rightarrow \Pi_{\infty}^{N}.\label{4.5.2}
\end{eqnarray}
This proves the first part of the theorem.

Define $r_{N}:(\mathbb{R}^{d})^{N} \to \mathcal{P}(\mathbb{R}^{d})$  as 
$$ r_{N}(x_{1},...x_{N}) = \frac{1}{N}\sum_{i=1}^{N}\delta_{x_{i}}, (x_1, \cdots x_N) \in (\R^d)^N.$$
Let $\nu_{n}^{N} = \Pi_{n}^{N}  \circ r_{N}^{-1}$ and $\nu_{\infty}^{N} = \Pi_{\infty}^{N}  \circ r_{N}^{-1}$. In order to prove that $\Pi^N_{\infty}$ is
$\mu_{\infty}$-chaotic, it suffices to argue that (cf. \cite{sznitman1991topics})
\begin{equation}
	\nu^N_{\infty} \to \delta_{\mu_{\infty}} \mbox { in } \mathcal{P}(\mathcal{P}(\mathbb{R}^{d})), \mbox{ as } N \to \infty.
\label{eq:eq727}	
\end{equation}
We first argue that as $n\to \infty$ 
\begin{equation}\nu_{n}^{N} \rightarrow \nu_{\infty}^{N}\quad\quad \text{in $\mathcal{P}(\mathcal{P}(\mathbb{R}^{d}))$}.\label{eq:eq728}\end{equation}
It suffices to show that $\langle F,\nu_{n}^{N}\rangle \to \langle F,\nu_{\infty}^{N}\rangle$ for any continuous and bounded function $F:\mathcal{P}(\mathbb{R}^{d}) \to \mathbb{R}$.
But this is immediate on observing that 
$$\langle F, \nu_n^N \rangle = \langle F\circ r_N , \Pi_n^N\rangle, \; \langle F, \nu_{\infty}^N \rangle = \langle F\circ r_N , \Pi_{\infty}^N\rangle,$$
the continuity of the map $r_N$ and the weak convergence of $\Pi_n^N$ to $\Pi^N_{\infty}$.
Next, for any 
$f \in BL_{1}(\mathcal{P}(\mathbb{R}^{d}))$
$$
|\langle f,\nu_{n}^{N} \rangle - \langle f,\delta_{\mu_{\infty}} \rangle| = | \frac{1}{n} \sum_{j=1}^n Ef(\mu_j^N) - f(\mu_{\infty})| 
\le \frac{1}{n} \sum_{j=1}^n E \mathcal{W}_{1}(\mu_{j}^{N},\mu_{\infty}).$$
Fix $\epsilon > 0$.  For every $N \in \mathbb{N}$ there exists $n_0(N) \in \N$ such that for all $n \ge n_0(N)$
$$
E \mathcal{W}_{1}(\mu_{n}^{N},\mu_{\infty}) \le \limsup_{n \to \infty} E \mathcal{W}_{1}(\mu_{n}^{N},\mu_{\infty}) + \epsilon.$$
Thus for all $n, N \in \mathbb{N}$
\begin{equation}
	|\langle f,\nu_{n}^{N} \rangle - \langle f,\delta_{\mu_{\infty}} \rangle| \le \frac{n_0(N)}{n} \max_{1 \le j \le n_0(N)} E \mathcal{W}_{1}(\mu_{j}^{N},\mu_{\infty})
	+ \limsup_{n \to \infty} E \mathcal{W}_{1}(\mu_{n}^{N},\mu_{\infty}) + \epsilon.
	\label{eq:eq740}
	\end{equation}
Finally
\begin{align*}
	\limsup_{N \to \infty} |\langle f,\nu_{\infty}^{N} \rangle - \langle f,\delta_{\mu_{\infty}} \rangle|
	= & 	\limsup_{N \to \infty} \lim_{n \to \infty}|\langle f,\nu_{n}^{N} \rangle - \langle f,\delta_{\mu_{\infty}} \rangle|  \\
	\le & \limsup_{N \to \infty} \limsup_{n \to \infty} E \mathcal{W}_{1}(\mu_{n}^{N},\mu_{\infty}) + \epsilon\\
	\le & \epsilon,
\end{align*}
where the first equality is from \eqref{eq:eq728}, the second uses \eqref{eq:eq740} and the third is a consequence of Corollary \ref{cor1}.  Since $\epsilon > 0$ is arbitrary, we have
\eqref{eq:eq727}	and the result follows.
\qed

\subsection{Proof of Theorem \ref{th5}}
\label{sec:secprfth5}
We will  first develop a concentration bound for $\mathcal{W}_{1}(\mu_{n}^{N},\Psi(\mu_{n-1}^{N}))$ for each fixed $n$ and then combine it with the estimate in
\eqref{4.3.2} in order to obtain the desired result.  The first step is carried out in the lemma below, the proof of which is given in Section \ref{sec:prooflem-4.20}.
\begin{lemma}
	\label{lem:fixn-4.10}
	Suppose Assumptions \ref{assumC} and \ref{assumD} hold. Then,
	there exist $a_1, a_2, a_3 \in (0, \infty)$ such that for all $\eps, R > 0$ and $n \in \mathbb{N}$,
	\begin{eqnarray*}
	P[\mathcal{W}_{1}(\mu_{n}^{N},\Psi(\mu_{n-1}^{N})) >\varepsilon] &\leq& a_3 \left(e^{-a_2 \frac{N \varepsilon^{2}}{R^{2}}} + \frac{R^{-\alpha}}{\varepsilon}\right)
	\end{eqnarray*}
	for all $N \ge \max\{1, a_1 (\frac{R}{\eps})^{d+2}\}$.  
\end{lemma}
   We now complete the proof of Theorem \ref{th5} using the lemma.
\subsubsection{Proof of Theorem \ref{th5}}
We will make use of \eqref{4.3.2}.  Recall that $\chi = e^{-\omega} + 2\delta \sigma$ and by our assumption $\chi \in (0,1)$.  
Let $\gamma = 2\sigma \gamma_0$.  Note that $\gamma < 1- e^{-\omega}$.
Then
\begin{eqnarray}
P[\mathcal{W}_{1}(\mu_{n}^{N},\mu_{n}) >\varepsilon] &\leq& P[\cup_{i=1}^{n}\{ \chi^{n- i}\mathcal{W}_{1}(\mu_{i}^{N},\Psi(\mu_{i -1}^{N})) >\gamma (1-\gamma)^{n-i}\varepsilon \} \cup \nonumber\\ &&\{\chi^{n}\mathcal{W}_{1}(\mu_{0}^{N},\mu_{0}) > \gamma (1-\gamma)^{n}\varepsilon\}]  \nonumber\\
&\leq& \sum_{i = 1}^{n} P[\mathcal{W}_{1}(\mu_{i}^{N},\Psi(\mu_{i -1}^{N})) >\gamma(\frac{1-\gamma}{\chi})^{n-i}\varepsilon] \nonumber\\&+& P[\mathcal{W}_{1}(\mu_{0}^{N},\mu_{0}) > \gamma (\frac{1-\gamma}{\chi})^{n}\varepsilon].\label{unif}
\end{eqnarray}
Let  $\beta = \gamma \eps$.  Note that $\vartheta = \frac{1-\gamma}{\chi}$ and
from our choice of $\delta$, $\vartheta > 1$.  Therefore
$$ N \ge a_1 (\frac{R}{\beta})^{d+2} \vee 1 \mbox{ implies } N \ge a_1 (\frac{R}{\beta \vartheta^n})^{d+2} \vee 1 \mbox{ for all } n \in \mathbb{N}_0.$$
Thus from Lemma \ref{lem:fixn-4.10}, for all $ N \ge a_1 (\frac{R}{\beta})^{d+1} \vee 1$ and $k=1, \cdots, n$
$$
P[\mathcal{W}_{1}(\mu_{k}^{N},\Psi(\mu_{k-1}^{N})) >\beta \vartheta^{n-k} ] \leq a_3 \left(e^{-a_2 \frac{N \beta^2\vartheta^{2(n-k)}}{R^{2}}} + \frac{R^{-\alpha}}{\beta \vartheta^{n-k}}\right).
$$
Using the above estimate in \eqref{unif}
\begin{align}
P[\mathcal{W}_{1}(\mu_{n}^{N},\mu_{n}) >\varepsilon] \le &
a_3 \sum_{i=0}^{n-1} \left(e^{-a_2 \frac{N \beta^2\vartheta^{2i}}{R^{2}}} + \frac{R^{-\alpha}}{\beta \vartheta^{i}}\right) + P[\mathcal{W}_{1}(\mu_{0}^{N},\mu_{0}) >\gamma\varepsilon] \nonumber	\\
\le & a_3 \sum_{i=0}^{\infty} e^{-a_2 \frac{N \beta^2\vartheta^{2i}}{R^{2}}} + \frac{a_3 R^{-\alpha}\vartheta}{\beta (\vartheta - 1)} + 
P[\mathcal{W}_{1}(\mu_{0}^{N},\mu_{0}) >\gamma\varepsilon].\label{eq:eq421}
\end{align}
Since $\vartheta >1$  we can find  $m_0 = m_{0}(\vartheta) \in \mathbb{N}$ such that
 $$\vartheta^{2i} \geq i \vartheta^{2} \quad\quad \forall i\ge m_{0}(\vartheta).$$ 
Thus 
\begin{eqnarray}
\sum_{i=0}^{\infty} e^{-a_2 \frac{N \beta^2\vartheta^{2i}}{R^{2}}} &=& \sum_{i=1}^{m_{0}(\vartheta)}e^{-a_2 \frac{N \beta^{2}\vartheta^{2i}}{R^{2}}} + \sum_{m_{0}(\vartheta)+1}^{\infty}e^{-a_2 \frac{N \beta^{2}\vartheta^{2i}}{R^{2}}}\nonumber\\
&\leq& m_{0}(\vartheta) e^{-a_2 \frac{N \beta^{2}\vartheta^{2}}{R^{2}}}+ \sum_{i=m_{0}(\vartheta)+1}^{\infty} (e^{-a_2 \frac{N \beta^{2}\vartheta^{2}}{R^{2}}})^{i}\nonumber\\ 
&\leq& [m_{0}(\vartheta) + \frac{1}{1-e^{-a_2 \frac{N \beta^{2}\vartheta^{2}}{R^{2}}}}]e^{-a_2 \frac{N \beta^{2}\vartheta^{2}}{R^{2}}}.\label{eq:eqab644}
\end{eqnarray}
Now for fixed $N \ge 1$ choose $R = \frac{\gamma \eps N^{1/d+2}}{a_1^{1/d+2}}$.  Then \eqref{eq:eqab644} holds for all such $N,R$.
Let $N_0 \ge 1$ be large enough so that for all $N \ge N_0$
\begin{equation}\label{eq:eq717}
1-e^{-a_2 a_1^{\frac{2}{d+2}}\vartheta^{2}N^{\frac{d}{d+2}}} > 1/2.\end{equation}
Then letting
$$ a_4 = a_3(m_0(\vartheta) + 2), \;\; a_5 = a a_1^{\frac{2}{d+2}}\vartheta^2, \;\; a_6 = \frac{\vartheta}{\vartheta-1} a_3 a_1^{\frac{\alpha}{d+2}} \gamma^{-(\alpha+1)},$$
we have for all $N \ge N_0$
$$
P[\mathcal{W}_{1}(\mu_{n}^{N},\mu_{n}) >\varepsilon] \le 
a_4 e^{-a_5N^{\frac{d}{d+2}}} + a_6 \eps^{-(\alpha+1)}N^{-\frac{\alpha}{d+2}} + P[\mathcal{W}_{1}(\mu_{0}^{N},\mu_{0}) >\gamma\varepsilon] .
$$
Choose $N_1 \ge N_0$ such that for all $N \ge N_1$, $N^{\frac{d}{d+2}} \ge \frac{2\alpha}{a_5 (d+2)}\log N$.
Also let $a_7 = \frac{2(1+\alpha)}{a_5}$.  Then for all $N \ge \max(N_1, (a_7 \log^+ \eps)^{(d+2)/d})$
$$
P[\mathcal{W}_{1}(\mu_{n}^{N},\mu_{n}) >\varepsilon] \le
(a_4 + a_6) \eps^{-(\alpha+1)} N^{-\frac{\alpha}{d+2}} + P[\mathcal{W}_{1}(\mu_{0}^{N},\mu_{0}) >\gamma\varepsilon].
$$
The result follows. \qed

\subsubsection{Proof of Lemma \ref{lem:fixn-4.10}.}
\label{sec:prooflem-4.20}
We now complete the proof of Lemma \ref{lem:fixn-4.10}.  The proof uses certain truncation ideas from \cite{bolley2007quantitative}.
Fix $\eps >0$.  For $\mu \in \mathcal{P}(\mathbb{R}^d)$, $R > 0$ and $\nu_0 \in \mathcal{P}(\mathbb{B}_R(0))$, where $\mathbb{B}_R(0) = \{x \in \mathbb{R}^d: |x| \le R\}$, define
$\mu_R \in \mathcal{P}(\mathbb{B}_R(0))$ as 
$$\mu_R(A) = \frac{\mu(A)}{\mu(\mathbb{B}_R(0))} 1_{\{\mu(\mathbb{B}_R(0)) \neq 0\}} + \nu_0(A)1_{\{\mu(\mathbb{B}_R(0)) = 0\}}, \;\;\; A \in \mathcal{B}(\mathbb{B}_R(0)).$$
For $N, n \in \mathbb{N}$ and $R > 0$, let  
$\Psi^{(R)}(\mu_{n-1}^{N}) := \frac{1}{N} \sum_{i=1}^{N} (\delta_{X_{n-1}^{i,N}}P^{\mu_{n-1}^{N}})_{R}$.

Let $\{Y_n^i\}_{i=1}^N$ be $\mathbb{B}_R(0)$ valued random variables which, conditionally on $\mathcal{F}_{n-1}^N$ are mutually independent and also independent of 
$\{X_n^{i,N}\}_{i=1}^N$, and 
$$
P(Y_n^i \in A \mid \mathcal{F}_{n-1}^N) = (\delta_{X_{n-1}^{i,N}} P^{\mu_{n-1}^{N}})_{R}(A), \;\; A \in \mathcal{B}(\mathbb{B}_R(0)).$$
Define
\[
Z_{n}^i = \left \{ 
\begin{array}{cc}
X_{n}^{i,N} & \text{ when }|X_{n}^{i,N}| \leq  R, \\ \ \\
Y_{n}^{i} & \text{ otherwise }.
\end{array} \right.
\]
It is easily checked that $P(Y_n^i \in A \mid \mathcal{F}_{n-1}^N) = P(Z_n^i \in A \mid \mathcal{F}_{n-1}^N)$ for all $A$ and conditionally on $\mathcal{F}_{n-1}^N$,
$\{Z_n^i\}_{i=1}^N$ are mutually independent.  Define $\mu_{n,R}^{N} :=  \frac{1}{N}\sum_{i=1}^{N} \delta_{Z_{n}^i}$.
Using triangle inequality we have
\begin{eqnarray}
\mathcal{W}_{1}(\mu_{n}^{N},\Psi(\mu_{n-1}^{N})) \leq \mathcal{W}_{1}(\Psi^{(R)}(\mu_{n-1}^{N}),\Psi(\mu_{n-1}^{N})) + \mathcal{W}_{1}(\Psi^{(R)}(\mu_{n-1}^{N}),\mu_{n,R}^{N}) + \mathcal{W}_{1}(\mu_{n}^{N},\mu_{n,R}^{N}).\nonumber\\
\label{4.6.1.2}
\end{eqnarray}
Consider first the middle term on the right side  of (\ref{4.6.1.2}). Recall $\mbox{Lip}^{0}_{1}(\mathbb{B}_{R}(0))=\{f \in \mbox{Lip}_{1}(\mathbb{B}_{R}(0)): f(0)=0\}.$ Then
\begin{eqnarray}
\mathcal{W}_{1}(\Psi^{(R)}(\mu_{n-1}^{N}),\mu_{n,R}^{N}) &=& \sup_{ f \in \mbox{Lip}^{0}_{1}(\mathbb{B}_R(0))} 
|\langle f,   \frac{1}{N} \sum_{i=1}^{N} (\delta_{X_{n-1}^{i,N}}P^{\mu_{n-1}^{N}})_{R} - \frac{1}{N} \sum_{i=1}^{N} \delta_{Z_n^i} \rangle|\nonumber\\
&=& \sup_{ f \in \mbox{Lip}^{0}_{1}(\mathbb{B}_R(0))} \left| \frac{1}{N} \sum_{i=1}^{N} \left(f(Z_n^i) - \langle f,(\delta_{X_{n-1}^{i,N}}P^{\mu_{n-1}^{N}})_{R} \rangle\right)   \right|\nonumber \\
&=& \sup_{ f \in \mbox{Lip}^{0}_{1}(\mathbb{B}_R(0))} |\frac{1}{N} \sum_{i=1}^{N} Z^{f}_{i,n}|\nonumber 
\end{eqnarray}
where $Z^{f}_{i,n} = f(Z_n^i) - \langle f,(\delta_{X_{n-1}^{i,N}}P^{\mu_{n-1}^{N}})_{R} \rangle.$  From Lemma \ref{app1}(a) there exists a finite subset $\mathcal{F}_{R,1}^{\frac{\varepsilon}{4}}(\mathbb{B}_R(0))$ of $\mbox{Lip}^{0}_{1}(\mathbb{B}_{R}(0))$ such that
\begin{eqnarray}\label{4.6.1.3}
\sup_{ f \in \mbox{Lip}^{0}_{1}(\mathbb{B}_R(0))} |\frac{1}{N} \sum_{i=1}^{N} Z^{f}_{i,n}|\leq \max_{f \in \mathcal{F}_{R,1}^{\frac{\varepsilon}{4}}(\mathbb{B}_R(0))}|\frac{1}{N} \sum_{i=1}^{N} Z^{f}_{i,n}|+ \frac{\varepsilon}{2}.
\end{eqnarray}
 Thus
\begin{eqnarray}
P[\mathcal{W}_{1}(\Psi^{(R)}(\mu_{n-1}^{N}),\mu_{n,R}^{N}) >\varepsilon] &\leq& E\left[\  P \left[\max_{f \in \mathcal{F}_{R,1}^{\frac{\varepsilon}{4}}(\mathbb{B}_R(0))}|\frac{1}{N} \sum_{i=1}^{N} Z^{f}_{i,n}| >\frac{\varepsilon}{2}\mid \mathcal{F}_{n-1}^{N}\right] \right] \nonumber\\
&\leq&  E\sum_{f \in \mathcal{F}_{R,1}^{\frac{\varepsilon}{4}}(\mathbb{B}_R(0))} P\left[ |\frac{1}{N} \sum_{i=1}^{N} Z^{f}_{i,n}| >\frac{\varepsilon}{2} \, \big\vert \, \mathcal{F}_{n-1}^{N} \right].
\end{eqnarray}
Since $f \in \mbox{Lip}^{0}_{1}(\mathbb{B}_R(0))$,  $|Z^{f}_{i,n}| \leq 2R$. So by  the Azuma - Hoeffding inequality, 
\begin{eqnarray}
P[\mathcal{W}_{1}(\Psi^{(R)}(\mu_{n-1}^{N}),\mu_{n,R}^{N}) >\epsilon] &\leq&  |\mathcal{F}_{R,1}^{\frac{\varepsilon}{4}}(\mathbb{B}_R(0))|
\max_{ f \in \mathcal{F}_{R,1}^{\frac{\varepsilon}{2}}(\mathbb{B}_R(0))}E\left( P[ |\frac{1}{N} \sum_{i=1}^{N} Z^{f}_{i,n}| >\frac{\varepsilon}{2}   \, \big\vert \, \mathcal{F}_{n-1}^{N}]\right)\nonumber\\
&\leq&  |\mathcal{F}_{R,1}^{\frac{\varepsilon}{4}}(\mathbb{B}_R(0))|2 e^{-\frac{N \varepsilon^{2}}{32 R^{2}}}\nonumber\\
&\leq&  2\ |\mathcal{F}_{R,1}^{\frac{\varepsilon}{4}}( [-R,R]^{d} )|e^{-\frac{N \varepsilon^{2}}{32R^{2}}}.
\end{eqnarray} 
From  Lemma \ref{app1}(b) we have the following estimate 
\begin{eqnarray}\label{mid}
P[\mathcal{W}_{1}(\Psi^{(R)}(\mu_{n-1}^{N}),\mu_{n,R}^{N}) >\varepsilon] &\leq& 
\max\left \{2, \frac{16R}{3\varepsilon}(2\sqrt{d}+1)3^{[\frac{8R}{\varepsilon}(\sqrt{d}+1)]^{d}  }\right\} e^{ -\frac{N \varepsilon^{2}}{32 R^{2}}}.
\end{eqnarray}
Thus there exist $k_1, k_2 \in (0,\infty)$ such that for all $n,N \in \mathbb{N}$, $R > 0, \eps > 0$
\begin{equation}
	\label{eq:eq1012}
P[\mathcal{W}_{1}(\Psi^{(R)}(\mu_{n-1}^{N}),\mu_{n,R}^{N}) >\varepsilon] \leq
k_2  [e^{k_1(R/\eps)^d}\vee 1] e^{ -\frac{N \varepsilon^{2}}{32 R^{2}}}.
\end{equation} 
For the first term in the right hand side of (\ref{4.6.1.2}) we make use of the observation that if for $i = 1, \cdots N$, $U^i, V^i$ are  $\mathbb{R}^d$ valued
random variables distributed according to $\lambda_U^i, \lambda_V^i$ respectively then
$$
\mathcal{W}_{1}(\frac{1}{N} \sum_{i=1}^N \lambda_U^i, \frac{1}{N} \sum_{i=1}^N \lambda_V^i) \le \frac{1}{N}\sum_{i=1}^N E|U^i-V^i|.$$
Thus
\begin{eqnarray*}
\mathcal{W}_{1}(\Psi^{(R)}(\mu_{n-1}^{N}),\Psi(\mu_{n-1}^{N})) &=&  \mathcal{W}_{1}(\frac{1}{N} \sum_{i=1}^{N} (\delta_{X_{n-1}^{i,N}}P^{\mu_{n-1}^{N}})_{R},\frac{1}{N} \sum_{i=1}^{N} \delta_{X_{n-1}^{i,N}} P^{\mu_{n-1}^{N}})\\
&\leq& \frac{1}{N} \sum_{i=1}^{N} E[|X_{n}^{i,N} - Z_n^i| \mid \mathcal{F}_{n-1}^{N}].
\end{eqnarray*}
Using the definition of $\{Z_n^i\}$ we see
\begin{eqnarray}
\frac{1}{N} \sum_{i=1}^{N} E[|X_{n}^{i,N} - Z_n^i| \mid \mathcal{F}_{n-1}^{N}]
&=& \frac{1}{N} \sum_{i=1}^{N} E[|X_{n}^{i,N} -Y_{n}^{i}|1_{|X_{n}^{i,N}| > R }  \mid \mathcal{F}_{n-1}^{N}]\nonumber\\
&\leq& \frac{2}{N} \sum_{i=1}^{N} E[|X_{n}^{i,N}|1_{|X_{n}^{i,N}| > R } \mid \mathcal{F}_{n-1}^{N}]\label{eq:eq1036}
 \end{eqnarray}
From \eqref{eq:eq258} we have that  $B_{n}^{\alpha} := E|X_{n}^{i,N}|^{1+\alpha}$ satisfies
$$B_{n}^{\alpha} \leq \kappa_1 E|X_0^{i,N}|^{1+\alpha} + \frac{\kappa_2}{1-\kappa_1} = B(\alpha).$$
Thus
\begin{eqnarray}
P[\mathcal{W}_{1}(\Psi^{(R)}(\mu_{n-1}^{N}), \Psi(\mu_{n-1}^{N})) > \varepsilon] &\leq& 
\frac{1}{\eps} E\frac{1}{N} \sum_{i=1}^{N} E[|X_{n}^{i,N} -Y_{n}^{i}|1_{|X_{n}^{i,N}| > R }  \mid \mathcal{F}_{n-1}^{N}]\nonumber\\
&\leq& \frac{2}{\varepsilon} E \{|X_{n}^{i,N}| 1_{|X_{n}^{i,N}| > R }\}\nonumber\\
&\leq& \frac{2 R^{-\alpha}}{\varepsilon} B_{n}^{\alpha} \leq \frac{2 R^{-\alpha}}{\varepsilon} B(\alpha).\label{first}
\end{eqnarray}
The third term in (\ref{4.6.1.2}) can be treated similarly.  Indeed, note that
\begin{eqnarray}
\mathcal{W}_{1}(\mu_{n}^{N},\mu_{n,R}^{N}) \leq \frac{1}{N} \sum_{i=1}^{N} |X_{n}^{i,N} - Z_n^i| =  \frac{1}{N} \sum_{i=1}^{N} |X_{n}^{i,N} -Y^{n}_{i}| 1_{|X_{n}^{i,N}| > R }. \nonumber
\end{eqnarray}
Thus using the bound for the right side of the first line in \eqref{first} we have that
\begin{align}
	P(\mathcal{W}_{1}(\mu_{n}^{N},\mu_{n,R}^{N}) > \eps) \leq & \frac{1}{\eps}E\frac{1}{N} \sum_{i=1}^{N} E[|X_{n}^{i,N} -Y^{n}_{i}|1_{|X_{n}^{i,N}| > R }  \mid \mathcal{F}_{n-1}^{N}]\nonumber\\
	\leq & \frac{2 R^{-\alpha}}{\varepsilon} B(\alpha).\label{eq:eq1051}
	\end{align}
	Using \eqref{eq:eq1012}, \eqref{first} and \eqref{eq:eq1051}  in \eqref{4.6.1.2}  we have 
\begin{eqnarray*}
P[\mathcal{W}_{1}(\mu_{n}^{N},\Psi(\mu_{n-1}^{N})) >\varepsilon] &\leq& 
k_2  [e^{k_1(3R/\eps)^d}\vee 1] e^{ -\frac{N \varepsilon^{2}}{288 R^{2}}}
+ \frac{12 R^{-\alpha}}{\varepsilon} B(\alpha).
\end{eqnarray*}
Letting $k_3 = 3^d\cdot 576 k_1$, $k_4 = 1/576$ and $k_5 = \max\{k_2, 12 B(\alpha)\}$, we have that
$$
P[\mathcal{W}_{1}(\mu_{n}^{N},\Psi(\mu_{n-1}^{N})) >\varepsilon] \leq
k_5  \left( e^{ -\frac{k_4 N \varepsilon^{2}}{ R^{2}}} 
+ \frac{ R^{-\alpha}}{\varepsilon} \right)
$$
for all $N \ge \max\{1, k_3 (\frac{R}{\eps})^{d+2}\}$.
This completes the proof of the lemma. \qed

\subsection{Proof of Theorem \ref{th6}}
We will proceed as in Section 
\ref{sec:secprfth5}  by first first giving a concentration bound for $\mathcal{W}_{1}(\mu_{n}^{N},\Psi(\mu_{n-1}^{N}))$ for each fixed $n$ and then combining
 it with 
\eqref{4.3.2} in order to obtain a uniform in $n$ estimate.  
We begin by observing that from Assumption \ref{assumG} it follows that there is a $\alpha_0 \in (0, \alpha]$ and $c_2 \in (0, \infty)$ such that 
for all $\alpha_1 \in [0, \alpha_0]$
\begin{equation}
	\label{eq:eq728b}
	\cle_1(\alpha_1) := \int e^{\alpha_1 D_1(z)} \theta(dz) \le e^{c_2\alpha_1}
\end{equation}
\begin{lemma}
	\label{lem:733}
	Suppose Assumptions \ref{assumC} and \ref{assumG} hold.  Let $\gamma_0$ be as in Theorem \ref{th6}.  Then for all
	$\delta \in [0,\min\{a_0-\gamma_0,\frac{1-e^{\omega}}{2M}\})$ and $\alpha_1 \in [0, \alpha_0]$
	$$\sup_{n\ge 0} \sup_{N \ge 1} E e^{\alpha_1|X_n^{1,N}|} < \infty .$$
\end{lemma}
{\bf Proof.}
Note that for $n \ge 1$
$$
|X_n^{i,N}| \le e^{-\omega} |X_{n-1}^{i,N}| + \delta M \left( |X_{n-1}^{i,N}| + \|\mu_{n-1}^N\|_1\right) + \delta D_1(\varepsilon^i_n).$$
Using Holder's inequality and taking expectations, for all $\alpha_1 \in [0, \alpha_0]$
\begin{align*}
	E e^{\alpha_1 |X_n^{i,N}|} & \le E \exp\left\{\alpha_1\left(e^{-\omega} |X_{n-1}^{i,N}| + \delta M \left( |X_{n-1}^{i,N}| + \|\mu_{n-1}^N\|_1\right) + \delta D_1(\varepsilon^i_n)\right)
	\right\}\\
	& \le \cle_1(\alpha_1 \delta) E \exp\left\{\alpha_1\left(e^{-\omega} |X_{n-1}^{i,N}| + \delta M \left( |X_{n-1}^{i,N}| + \|\mu_{n-1}^N\|_1\right)\right)
	\right\}\\
	&\le \cle_1(\alpha_1 \delta) E \exp\left\{\alpha_1 (e^{-\omega}+ 2\delta M) |X_{n-1}^{i,N}|\right\},
\end{align*}
where the last inequality is from Jensen's inequality.
Thus for all $\alpha_1 \in [0, \alpha_0]$
$$f_n(\alpha_1) := E \exp\{\alpha_1 |X_n^{i,N}|\}  \le \cle_1(\alpha_1 \delta) f_{n-1}(\alpha_1 \kappa_1),$$
where by our assumption $\kappa_1 = e^{-\omega} + 2 \delta M \in (0,1)$.  Iterating the above inequality we have for all $n \ge 1$
$$f_n(\alpha_1) \le f_0(\alpha_1) \prod_{j=0}^{n-1} \cle_1(\alpha_1 \delta \kappa_1^j)
\le f_0(\alpha_1) e^{c_2\alpha_1 \delta \sum_{j=0}^{n-1}\kappa_1^j} \le f_0(\alpha_1) e^{c_2\alpha_1 \delta/(1-\kappa_1)}$$
where the second inequality is a consequence of \eqref{eq:eq728b}.
The result follows. \qed

The following lemma is proved
in a manner similar to Lemma \ref{lem:fixn-4.10} so only a sketch is provided.
\begin{lemma}\label{lem:lem352}
	There exist $\ti a_1, \ti a_2, \ti a_3 \in (0, \infty)$ and, for each $\alpha_1 \in [0, \alpha_0)$,  $\ti B(\alpha_1) \in [0, \infty)$ such that for all $\eps, R > 0$ and $n \in \mathbb{N}$, 
	\begin{eqnarray*}
	P[\mathcal{W}_{1}(\mu_{n}^{N},\Psi(\mu_{n-1}^{N})) >\varepsilon] &\leq& \ti a_3 \left(e^{-\ti a_2 \frac{N \varepsilon^{2}}{R^{2}}} + \ti B(\alpha_1)\frac{e^{-\alpha_1 R}}{\varepsilon}\right)
	\end{eqnarray*}
	for all $N \ge \max\{1, \ti a_1 (\frac{R}{\eps})^{d+2}\}$.	
\end{lemma}
{\bf Proof.}  From Lemma \ref{lem:733} we have that for $\alpha_1 \in [0, \alpha_0]$
\begin{equation}
	\label{eq:eq553}
	\sup_{n \ge 0} \sup_{N \ge 1} \max_{\{1 \le i \le N\}} E e^{\alpha |X_n^{i,N}|} < \infty .
\end{equation}
Next, as in the proof of Lemma \ref{lem:fixn-4.10}, we will use \eqref{4.6.1.2}.  For the middle term on the right side of \eqref{4.6.1.2} we use the same bound as in \eqref{eq:eq1012}.
Now consider the first term in \eqref{4.6.1.2}.  From \eqref{eq:eq1036} we have that
\begin{eqnarray}
P[\mathcal{W}_{1}(\Psi(\mu_{n-1}^{N}),\Psi(\mu_{n-1}^{N})_{R})>\varepsilon] &\leq& \frac{2}{\varepsilon}E\left(|X_{n}^{1,N}| 1_{|X_{n}^{1,N}|>R}\right).\label{eq:eq603}
\end{eqnarray}
From \eqref{eq:eq553} it follows that for every $\alpha_1 \in [0, \alpha_0)$
\begin{equation*}
	\sup_{n \ge 0} \sup_{N \ge 1} \max_{\{1 \le i \le N\}} E \left(|X_n^{i,N}|e^{\alpha_1 |X_n^{i,N}|}\right) = \ti B(\alpha_1) < \infty .
\end{equation*}
Applying Markov's inequality we now have for $\alpha_1 \in [0, \alpha_0)$
\begin{eqnarray}
P[\mathcal{W}_{1}(\Psi(\mu_{n-1}^{N}),\Psi(\mu_{n-1}^{N})_{R})>\varepsilon] &\leq& \frac{2}{\varepsilon} e^{-\alpha_1 R}\ti B(\alpha_1).\label{eq:eq613}
\end{eqnarray}
The third term in \eqref{4.6.1.2} is bounded similarly.  Indeed, as in \eqref{eq:eq1051} we get for $\alpha_1 \in [0, \alpha)$
\begin{equation}
	P[\mathcal{W}_{1}(\mu_{n}^{N},\mu_{n,R}^{N}) > \eps] \leq \frac{2}{\varepsilon} e^{-\alpha_1 R}\ti B(\alpha_1).\label{eq:eq616}
\end{equation}
Using \eqref{eq:eq1012}, \eqref{eq:eq613} and \eqref{eq:eq616} in \eqref{4.6.1.2} we now have for $\alpha_1 \in [0, \alpha_0)$
\begin{eqnarray*}
P[\mathcal{W}_{1}(\mu_{n}^{N},\Psi(\mu_{n-1}^{N})) >\varepsilon] &\leq& 
k_2  [e^{k_1(3R/\eps)^d}\vee 1] e^{ -\frac{N \varepsilon^{2}}{288 R^{2}}}
+ \frac{12 e^{-\alpha_1 R}}{\varepsilon} \ti B(\alpha_1).
\end{eqnarray*}
Thus with $k_3, k_4$ as in the proof of Lemma \ref{lem:fixn-4.10} and $k_5 = \max\{k_2, 12\}$ we have
$$
P[\mathcal{W}_{1}(\mu_{n}^{N},\Psi(\mu_{n-1}^{N})) >\varepsilon] \leq
k_5  \left( e^{ -\frac{k_4 N \varepsilon^{2}}{ R^{2}}} 
+ \ti B(\alpha_1)\frac{ e^{-\alpha_1 R}}{\varepsilon} \right)
$$
for all $N \ge \max\{1, k_3 (\frac{R}{\eps})^{d+2}\}$.
The result follows. \qed

We now complete the proof of Theorem \ref{th6}.
\subsubsection{Proof of Theorem \ref{th6}.}
Fix $\alpha_1 \in [0, \alpha_0)$.
Following the steps  in  the proof of \eqref{eq:eq421}, with $\vartheta, \beta$ as in Theorem \ref{th5}, we have from Lemma \ref{lem:lem352}, for all $ N \ge \ti a_1 (\frac{R}{\beta})^{d+2} \vee 1$ and $k=1, \cdots, n$
\begin{align}
P[\mathcal{W}_{1}(\mu_{n}^{N},\mu_{n}) >\varepsilon] \le
 \ti a_3 \sum_{i=0}^{\infty} e^{-\ti a_2 \frac{N \beta^2\vartheta^{2i}}{R^{2}}} + \frac{\ti a_3 \ti B(\alpha_1) e^{-\alpha_1 R}\vartheta}{\beta (\vartheta - 1)} + 
P[\mathcal{W}_{1}(\mu_{0}^{N},\mu_{0}) >\gamma\vartheta^n\varepsilon].\label{eq:eq635}
\end{align}
As before for fixed $N \ge 1$ choose $R = \frac{\gamma \eps N^{1/d+2}}{\ti a_1^{1/d+2}}$.  Then \eqref{eq:eqab644} holds for all such $N,R$ with $a_2$ replaced by $\ti a_2$.
Let $N_0 \ge 1$ be large enough so that for all $N \ge N_0$, \eqref{eq:eq717} holds with $(a_1, a_2)$ replaced by $(\ti a_1, \ti a_2)$.
%
%
%
Then letting
$$ \ti a_4 = \ti a_3(m_0(\vartheta) + 2), \;\; \ti a_5 = \ti a_2 \ti a_1^{\frac{2}{d+2}}\vartheta^2, \;\; \ti a_6 = \frac{\ti a_3 \vartheta}{\gamma(\vartheta-1)} \ti B(\alpha_1),\;\;
\ti a_7 = \frac{\alpha_1 \gamma}{\ti a_1^{1/d+2}},$$
we have for all $N \ge N_0$
$$
P[\mathcal{W}_{1}(\mu_{n}^{N},\mu_{n}) >\varepsilon] \le 
\ti a_4 e^{-\ti a_5N^{\frac{d}{d+2}}} + \ti a_6 \eps^{-1}\exp (-\ti a_7 \e N^{\frac{1}{d+2}}) + P[\mathcal{W}_{1}(\mu_{0}^{N},\mu_{0}) >\gamma\vartheta^n\varepsilon] .
$$
Note that $\eps^{-1}\exp (-\ti a_7 \e N^{\frac{1}{d+2}}) < \exp (-\frac{\ti a_7}{2} \e N^{\frac{1}{d+2}}) $ if $N > \left(\frac{2}{\ti a_7}\right)^{d+2}(\frac{1}{\e} \log^+ \frac{1}{\e})^{d+2}$. 

Consider now the case $d > 1$.  Then, taking $L_1 = \max\{(\frac{2}{\ti a_7})^{d+2}, N_0\}$, $L_2 = \ti a_4 + \ti a_6$, 
$L_3 = \min \{\ti a_5, \ti a_7/2\}$, we have
for all $N \ge L_1 \max\{ (\frac{1}{\e} \log^+ \frac{1}{\e})^{d+2}, \e^{(d+2)/(d-1)}\}$
$$
P[\mathcal{W}_{1}(\mu_{n}^{N},\mu_{n}) >\varepsilon] \le L_2 e^{-L_3 \e N^{1/d+2}}+ P[\mathcal{W}_{1}(\mu_{0}^{N},\mu_{0})>\gamma\vartheta^n\varepsilon].$$
This proves the theorem for the case $d> 1$.
Finally for $d=1$, with the same choice of $L_1, L_2, L_3$, we have for all 
$N \ge L_1 \max\{ (\frac{1}{\e} \log^+ \frac{1}{\e})^{d+2}, 1\}$
$$
P[\mathcal{W}_{1}(\mu_{n}^{N},\mu_{n}) >\varepsilon] \le L_2 e^{-L_3 (\e \wedge 1) N^{1/d+2}}.$$
The result follows.
 \qed

\subsection{Proof of Theorem \ref{thm6}.}
\label{sec:secprfthm6}
In order to prove the theorem we will introduce an auxiliary sequence $\{Y_n^{i,N}, i=1, \cdots N\}_{n\ge 0}$ such that for each $n$, $\{Y_n^{i,N}\}_{i=1}^N$ are i.i.d.  We will then employ results from \cite{boissard2011simple} and \cite{bolley2005weighted} in order to give a uniform (in $k$) concentration bound for
$\mathcal{W}_{1}(\eta_{k}^{N},\mu_{k})$, where $\eta_{k}^{N}$ is the empirical measure $\frac{1}{N} \sum_{k=1}^N \delta_{Y_k^{i,N}}$.
Finally we will obtain the desired concentration estimate on $\mathcal{W}_{1}(\mu_{k}^{N},\mu_{k})$ by making use of Lemma \ref{lem:lemgron} below.  We begin by introducing our auxiliary system.
\subsubsection{An Auxiliary System.}
\label{sec:secaux}
Consider the collection of $\RR^d$ valued random variables $\{Y_{n}^{i,N},  i=1,...,N\}_{n\geq 0}$ defined as follows.
\begin{eqnarray}
Y_{n+1}^{i,N} &=& A Y_{n}^{i,N} + \delta f(Y_{n}^{i,N},\mu_{n},\epsilon_{n+1}^{i}), \;\; n \ge 0\nonumber\\
Y_{0}^{i,N} &=& X_{0}^{i,N}.\label{add}
\end{eqnarray}
Note that for each $n$, $\{Y_n^{i,N}\}_{i=1}^N$ are i.i.d.
 In fact, since  $\mathcal{L}(\{X_{0}^{i,N}\}_{i=1,...N}) = \mu_{0}^{\otimes N}$, we have $\mathcal{L}(\{Y_{n}^{i,N}\}_{i=1,...N}) = \mu_{n}^{\otimes N}.$
Let $\eta_{n}^{N} := \frac{1}{N}\sum_{i=1}^{N} \delta_{Y_{n}^{i,N}}.$
The following lemma will give a useful relation between $\mathcal{W}_{1}(\eta_{n}^{N},\mu_{n})$ and $\mathcal{W}_{1}(\mu_{n}^{N},\mu_{n}).$
\begin{lemma}
	\label{lem:lemgron}
	Suppose Assumptions  \ref{assumC} and \ref{assumG} hold.
	Let $\chi_2 = e^{-\omega} + 2\delta M$. Then for every $n \ge 0$ and $N \ge 1$
	\begin{eqnarray}
	\mathcal{W}_{1}(\mu_{n+1}^{N},\mu_{n+1}) &\leq& \mathcal{W}_{1}(\eta_{n+1}^{N},\mu_{n+1}) +\delta M\sum_{k=0}^{n} \chi_{2}^{n-k} \mathcal{W}_{1}(\eta_{k}^{N},\mu_{k}).\label{Gron}
	\end{eqnarray}
\end{lemma}
\textbf{Proof:}
Since by Assumption \ref{assumG}  $D(\epsilon)\leq M$, we have for each $i=1,...,N$
\begin{eqnarray*}
|X_{n+1}^{i,N} - Y_{n+1}^{i,N}| &\leq& e^{-\omega} |X_{n}^{i,N} - Y_{n}^{i,N}| + \delta M \{|X_{n}^{i,N} - Y_{n}^{i,N}| +\mathcal{W}_{1}(\mu_{n}^{N},\mu_{n})\}\\
&=&(e^{-\omega} + \delta M)|X_{n}^{i,N} - Y_{n}^{i,N}| +  \delta M  \mathcal{W}_{1}(\mu_{n}^{N},\mu_{n})\\
\end{eqnarray*}
Thus
\begin{equation}
|X_{n+1}^{i,N} - Y_{n+1}^{i,N}| \le \delta M \sum_{k=0}^{n}(e^{-\omega}+\delta M )^{n-k} \mathcal{W}_{1}(\mu_{k}^{N},\mu_{k}). \label{eq:eq612}
\end{equation} 
Now note that 
\begin{eqnarray}
\mathcal{W}_{1}(\eta_{n+1}^{N},\mu_{n+1}^{N}) &\leq& \frac{1}{N} \sum_{i=1}^{N}|X_{n+1}^{i,N} - Y_{n+1}^{i,N}|\leq \delta M \sum_{k=0}^{n}(e^{-\omega}+\delta M )^{n-k}\mathcal{W}_{1}(\mu_{k}^{N},\mu_{k}).\nonumber 
\end{eqnarray}
Using triangle inequality
\begin{eqnarray}
\mathcal{W}_{1}(\eta_{n+1}^{N},\mu_{n+1}^{N}) &\leq& \delta M \sum_{k=0}^{n}(e^{-\omega}+\delta M )^{n-k}\mathcal{W}_{1}(\eta_{k}^{N},\mu_{k}^{N})+\delta M \sum_{k=0}^{n}(e^{-\omega}+\delta M )^{n-k}\mathcal{W}_{1}(\eta_{k}^{N},\mu_{k}).\nonumber 
\end{eqnarray}
Applying Lemma \ref{app4} with
$$a_n = \chi_{1}^{-n} \mathcal{W}_{1}(\eta_{n}^{N},\mu_{n}^{N}), \; b_n = \frac{\delta M}{\chi_{1}} \sum_{k=0}^{n-1} \chi_{1}^{-k} \mathcal{W}_{1}(\eta_{k}^{N},\mu_{k}), \;
c_n = \frac{\delta M}{\chi_{1}}, \; n \ge 0$$
where 
$\chi_{1}:= e^{-\omega}+\delta M$, we have 
\begin{eqnarray}
\chi_{1}^{-(n+1)}\mathcal{W}_{1}(\eta_{n+1}^{N},\mu_{n+1}^{N}) &\leq& b_{n+1}+  \sum_{k=0}^{n}(\frac{\delta M}{\chi_{1}})^{2} \sum_{i=0}^{k-1} \chi_{1}^{-i}\mathcal{W}_{1}(\eta_{i}^{N},\mu_{i}) \left(1+\frac{\delta M}{\chi_{1}}\right)^{n-k}\nonumber\\
&=& b_{n+1}+ \sum_{i=0}^{n}\sum_{k=i+1}^{n}(\frac{\delta M}{\chi_{1}})^{2}(1+\frac{\delta M}{\chi_{1}})^{n-k}\chi_{1}^{-i}\mathcal{W}_{1}(\eta_{i}^{N},\mu_{i}) \nonumber\\
&=& b_{n+1}+\sum_{i=0}^{n} (\frac{\delta M}{\chi_{1}})^{2}\chi_{1}^{-i}.\mathcal{W}_{1}(\eta_{i}^{N},\mu_{i})\sum_{m=0}^{n-i-1} (1+\frac{\delta M}{\chi_{1}})^{m}\nonumber\\
&=& b_{n+1}+\sum_{i=0}^{n} (\frac{\delta M}{\chi_{1}})\chi_{1}^{-i}\mathcal{W}_{1}(\eta_{i}^{N},\mu_{i})[(1+\frac{\delta M}{\chi_{1}})^{n-i} -1]\label{Gron1}. 
\end{eqnarray}
Simplifying (\ref{Gron1}) one gets 
\begin{eqnarray}
\mathcal{W}_{1}(\eta_{n+1}^{N},\mu_{n+1}^{N}) &\leq& \delta M \sum_{k=0}^{n} \chi_{1}^{n-k}\mathcal{W}_{1}(\eta_{k}^{N},\mu_{k}) +\sum_{k=0}^{n}\delta M \chi_{1}^{n-k} \mathcal{W}_{1}(\eta_{k}^{N},\mu_{k})[(1+\frac{\delta M}{\chi_{1}})^{n-k} - 1]\nonumber\\
&=& \delta M\sum_{k=0}^{n} (\chi_{1}+\delta M)^{n-k}\mathcal{W}_{1}(\eta_{k}^{N},\mu_{k}).\nonumber
\end{eqnarray}
The result now follows by an application of triangle inequality.
\qed
\subsubsection{Transportation inequalities.}
\label{sec:sectransp}
Proof of Theorem \ref{thm6} is based on certain results from \cite{boissard2011simple} and \cite{bolley2005weighted} which we summarize in this section.
We begin with the definition of a `transportation inequality'.
Recall that given $\mu, \nu \in \mathcal{P}(\RR^d)$, the relative entropy of $\nu$ with respect to $\mu$ is defined as
\[
R(\nu || \mu) = \left \{ 
\begin{array}{cc}
\int_{\RR^d} \left(\log \frac{d\nu}{d\mu}\right) d\nu  & \text{ when } \nu \ll \mu, \\ \ \\
\infty & \text{ otherwise }.
\end{array} \right.
\]
Define $\ell: [0, \infty) \to [0, \infty)$ as $\ell(x) = x \log x - x + 1$.
\begin{definition}
	\label{def:transport}
	Let $\alpha: [0, \infty) \to \mathbb{R}$ be a convex, increasing left continuous function such that $\alpha(0) = 0$.
	We say $\nu \in \mathcal{P}(\RR^d)$ satisfies a $\alpha(\mathcal{T})$ inequality if for all $\tilde \nu \in \mathcal{P}(\RR^d)$
	$$\alpha (\mathcal{W}_1(\nu, \tilde \nu)) \le R(\tilde \nu || \nu).$$
\end{definition}
The following result is established in \cite{boissard2011simple}.
\begin{Theorem}
	\label{thm:thmboiss}
	Suppose that $\nu \in \mathcal{P}(\RR^d)$ satisfies a $\alpha(\mathcal{T})$ inequality and suppose that there is  $\zeta > 0$ such that
	$\int_{\RR^d} e^{\zeta |x|} \nu(dx) \le 2.$
	Let $L_N = \frac{1}{N} \sum_{i=1}^N \delta_{Z_i}$ where $Z_i$ are i.i.d. with common distribution $\nu$.
	Then for $t > 0$
	$$
	P(\mathcal{W}_1(L_N, \nu) \ge t) \le \exp \left\{-N \alpha(\frac{t}{2} - \Gamma(\mathcal{C}_t, N)) \right\},$$
	where
	\begin{equation}
		\Gamma(\mathcal{C}_t, N)) = \inf_{\lambda > 0} \left \{\frac{1}{\lambda} \log \mathcal{C}_t + N \alpha^*(\frac{\lambda}{N}) \right \},
		\label{eq:eq548}
	\end{equation}
	$\alpha^*: \RR \to [0, \infty)$ is defined as
	$$ \alpha^*(s) = \sup_{t \ge 0} \{st - \alpha(t) \} 1_{[0, \infty)}(s), \; s \in \R,$$
	\begin{equation}
		\label{eq:eq549}
		\mathcal{C}_t = 2 \left ( 1 + \psi(\frac{32}{\zeta t})\right) 2^{c_d (\psi(\frac{32}{\zeta t}))^d},
	\end{equation}
	$\psi(x) = x\log(2\ell(x))$, $x \ge 0$ and $c_d$ is a positive scaler depending only on $d$. 
\end{Theorem}
The following result is from \cite{bolley2005weighted}.
\begin{Theorem}
	\label{th:thboll}
	Let $\nu \in \mathcal{P}(\mathbb{R}^d)$. Suppose that $\int_{\RR^d} e^{\alpha_0 |x|} d\nu(x) < \infty$ for some $\alpha_0 > 0$.  Then $\nu$ satisfies $\alpha(\mathcal{T})$ inequality with
	\begin{equation}\label{eq:eq626}
		\alpha(t) = \left ( \sqrt{\frac{t}{C} + \frac{1}{4}} - \frac{1}{2}\right)^2, \; t \ge 0
	\end{equation}
	for any 
	$$C > 2  \inf_{x_0 \in \RR^d, \tilde \alpha > 0} \frac{1}{\tilde \alpha}\left ( \frac{3}{2} + \log \int_{\RR^d} e^{\tilde \alpha |x-x_0|} d \nu(x) \right).$$
\end{Theorem}
\subsubsection{Exponential Integrability.}
\label{sec:expinteg}
Transportation inequalities presented in Section \ref{sec:sectransp} require exponential integrability of the underlying measure.  In this section we show that under
Assumption \ref{assumG} the desired integrability properties hold.  
\begin{lemma}
	\label{lem:expinteg}
	Suppose that Assumption \ref{assumG} holds and that $\delta \in (0, \frac{1-e^{-\omega}}{2M})$.  Then
	$\kappa_1 = (e^{-\omega} + 2\delta M) \in (0,1)$ and for all $\alpha_1 \in [0, \alpha_0]$
		$$\sup_{n \ge 0} \int_{\RR^d} e^{\alpha_1 |x|} \mu_n(dx) \le (\int e^{\alpha_{1} |x|} \mu_{0}(dx) ) \exp\left\{\frac{c_2 \alpha_{1}}{1-\kappa_{1}}\right\}.$$
\end{lemma}
{\bf Proof.}
The property that $\kappa_1 \in (0,1)$ is an immediate consequence of assumptions on $\delta$.
Let
 $f_{n}(\alpha_1) := \int e^{\alpha_1 |x|}\mu_{n}(dx)$. 

From (\ref{eq1.1.1}) and the condition $D(\epsilon)\leq M$ we have 
\begin{eqnarray}
|X_{n+1}|\leq e^{-\omega}|X_{n}|+\delta M(|X_{n}| +||\mu_{n}||_{1})+\delta D_{1}(\epsilon_{n+1}).
\end{eqnarray}
Using Holder's inequality  and taking exponentials we get 
\begin{eqnarray}
f_{n+1}(\alpha_1)= E e^{\alpha_1 |X_{n+1}|}&\leq& Ee^{\alpha_1[e^{-\omega}|X_{n}|+\delta M(|X_{n}| +||\mu_{n}||_{1})+\delta D_{1}(\epsilon_{n+1})]}\nonumber\\
&=& \mathcal{E}_{1}(\alpha_1 \delta)  Ee^{\alpha_1 (e^{-\omega} +\delta M)|X_{n}| + \delta M||\mu_{n}||_{1}}. \label{genu}
\end{eqnarray}
From Jensen's inequality we have $||\mu_{n}||_{1}\leq E|X_{n}|$. Applying Jensen's inequality again to the  function 
$x \mapsto \exp\{\alpha_1 \delta M  x\}$
we have 
\begin{eqnarray}
f_{n+1}(\alpha_1)&\leq& \mathcal{E}_{1}(\alpha_1 \delta)  E[e^{\alpha_1 (e^{-\omega} +\delta M)|X_{n}| }]E[e^{\alpha_1 \delta M|X_{n}|}]. \label{genu1}
\end{eqnarray}
Note that for any two non-decreasing, non-negative functions $f,g$ on $\R$ and any $\pi \in \PP(\R)$,
  $$\int f(x)g(x) \mu(dx) \geq \int f(x) \pi(dx) \int g(y) \pi(dy).$$  Using this inequality in the above display  yields the following recursion
\begin{eqnarray*}
f_{n+1}(\alpha_1) \leq  \mathcal{E}_{1}(\alpha_1 \delta)  E[e^{\alpha_1 (e^{-\omega} +2\delta M)|X_{n}| }]
= \mathcal{E}_{1}(\alpha_1 \delta)  f_{n}(\alpha_1 \kappa_1).
\end{eqnarray*}
Iterating the above inequality we have, for all $n \ge 0$,
$$f_{n+1}(\alpha_1) \le f_0(\alpha_1)\prod_{j=0}^n \mathcal{E}_1(\alpha_1 \delta \kappa_1^j).$$
  Thus using \eqref{eq:eq728b} we see
\begin{align*}
	f_{n+1}(\alpha_1) \le  f_0(\alpha_1)\prod_{j=0}^n \exp \{ c_2 (\alpha_1 \delta \kappa_1^j)\}
	\le  f_0(\alpha_1) \exp \{ c_2 \alpha_1 \sum_{j=0}^{\infty} \kappa_1^j\}.
\end{align*}
The result follows.
\qed
\subsection{Uniform Concentration Bounds for $\{\eta_n^N\}$.}
In this section we will give, using results of Sections \ref{sec:sectransp} and \ref{sec:expinteg}, uniform concentration bounds for $\{\eta_n^N\}_{n \ge 1}$ as $N \to \infty$.
\begin{lemma}
	\label{cor:cor621}
	Suppose that Assumption \ref{assumG} holds and $\delta \in (0, \frac{1-e^{-\omega}}{2M})$.
	Then the following hold.
	\begin{description}
	\item{(1)} There exists a $\zeta_0 \in (0, \infty)$ such that 
	\begin{equation}
		\sup_{n \in \mathbb{N}_0} \int_{\RR^d} e^{\zeta_0|x|} \mu_n(dx) \le 2 \label{eq:eq703b}
	\end{equation}
	and for all $n \in \mathbb{N}_0$, $\mu_n$ satisfies a $\alpha(\mathcal{T})$ inequality with $\alpha$ as in \eqref{eq:eq626}
	and with
	$$C \ge C_0= 2 \sqrt{2}  \frac{1}{\zeta_0}\left ( \frac{3}{2} + \log 2 \right).$$
	\item{(2)} For all $t > 0$ and $n \in \mathbb{N}_0$
	$$
	P(\mathcal{W}_1(\eta_n^N, \mu_n) \ge t) \le \exp \left\{-N \alpha_0(\frac{t}{2} - \Gamma_0(\mathcal{C}^0_t, N)) \right\},$$
	where $\alpha_0$ is defined by the right side  in \eqref{eq:eq626} with $C$ replaced with $C_0$, $\Gamma_0$ is defined by the right side of 
	 \eqref{eq:eq548}
	with $\alpha^*$ replaced by $\alpha_0^*$ and $\mathcal{C}^0_t$ is as in \eqref{eq:eq549} with $\zeta$ replaced with $\zeta_0$.\\
	\item{(3)} There exist $N_1 \in \mathbb{N}$ and $L_1 \in (0, \infty)$ such that
for all $t \in [\frac{C_0}{2}, \infty)$, $n \in \N_0$ and $N \ge N_1$
	$$P(\mathcal{W}_1(\eta_n^N, \mu_n) \ge t) \le  \exp (-L_1 N t).$$
\item{(4)} There exist $L_2, L_3 \in (0, \infty)$ such that for all $t \in (0, \frac{C_0}{2}]$ and all $N \ge L_3 \frac{\log \mathcal{C}_t}{t^2}$.
	$$P(\mathcal{W}_1(\eta_n^N, \mu_n) \ge t) \le \exp (-L_2 Nt^2).$$
\end{description}
\end{lemma}
{\bf Proof.}
(1)   Suppose that the statement in \eqref{eq:eq703b} fails to hold for any $\zeta_0 > 0$.  Then there exist  sequences $n_k \uparrow \infty$ and $\zeta_k \downarrow 0$ such that
\begin{equation}
\int_{\RR^d} e^{\zeta_{k}|x|} \mu_{n_k}(dx) > 2. \label{eq:eq706}	
\end{equation}
From Lemma \ref{lem:expinteg} it follows that $\{\mu_{n_k}, k \ge 1\}$ is tight.  Suppose along a further subsequence $\mu_{n_k}$
 converges to some measure $\mu_0$.  Then sending $k \to \infty$ along this subsequence in \eqref{eq:eq706}	and using Lemma \ref{lem:expinteg} once again we arrive at a contradiction.
This proves the first statement in (1).  The second statement in (1) is an immediate consequence of Theorem \ref{th:thboll}.\\ 

(2) This is immediate from part (1) and Theorem \ref{thm:thmboiss}.\\ 

(3) It is easy to check that for all $t > 0$, $N \in \mathbb{N}$ (see proof of Corollary 2.5 in \cite{boissard2011simple})
$$\Gamma_0(\mathcal{C}^0_t, N) \le \frac{C_0}{ \left(1+ \frac{N}{\log\mathcal{C}^0_t}\right)^{1/2}-1}.$$
Thus recalling the expression for $\mathcal{C}^0_t$ in \eqref{eq:eq549} we see that $\lim_{N \to \infty} \sup_{t \ge C_0/2} \Gamma_0(\mathcal{C}^0_t, N) =0$.
Choose $N_1 \in \mathbb{N}$ such that for all $N \ge N_1$ and $t \ge C_0/2$
$$ \Gamma_0(\mathcal{C}^0_t, N) \le \frac{C_0}{8} \le \frac{t}{4}.$$
Then for all $N \ge N_1$ and $t \ge C_0/2$
\begin{align*}
	\alpha_0\left(\frac{t}{2} - \Gamma_0(\mathcal{C}^0_t, N)\right) \ge & \frac{1}{4} \left ( (1+ \frac{t}{C_0})^{1/2} -1 \right)^2
	\ge \frac{1}{16} \frac{(t/C_0)^2}{(1 + t/C_0)} \ge \frac{t}{48C_0},
\end{align*}
where the second inequality follows on using the inequality
\begin{equation} \label{eq: eq446}\sqrt{1+x} - 1 \ge \frac{x}{2\sqrt{1+x}}, \; x \ge 0.\end{equation}
Combining this with (2) completes the proof of (3).

(4)  From the proof of Corollary 2.5 of \cite{boissard2011simple} it follows that for $t \le C_0/2$
\begin{equation}\label{eq:eq805}
P(\mathcal{W}_1(\eta_n^N, \mu_n) \ge t)
\le A(N, t) \exp \left ( - B_1 N t^2\right )
\end{equation}
where $A(N,t) = \exp \left ( \frac{N B_2}{ \left ((1 + N/ \log \mathcal{C}^0_t)^{1/2} -1\right)^2}\right)$, $B_1 = (\sqrt{2}-1)^2/(2C_1^2)$ and
$B_2 = 4(\sqrt{2}-1)^2$.

From \eqref{eq: eq446} note that if $N > \log \mathcal{C}^0_t$
$$\left ((1 + N/ \log \mathcal{C}^0_t)^{1/2} -1\right)^2 \ge \frac{N}{8 \log \mathcal{C}^0_t}.$$
Thus
for all such $N,t$, $A(N,t) \le \exp (8B_2 \log \mathcal{C}^0_t)$.  Thus if additionally $N \ge \frac{16B_2}{B_1} \frac{\log \mathcal{C}^0_t}{t^2}$, the right side
of \eqref{eq:eq805} is bounded above by $\exp (-B_1 Nt^2/2)$.  The result follows. \qed

%
%

\subsubsection{Proof of Theorem \ref{thm6}.}
\label{sec:subsecthm6prf}
In this section we complete the proof of Theorem \ref{thm6}.  

Fix $\gamma \in (0, 1 - e^{-w})$.  From \eqref{Gron}, for any $\e > 0$,
\begin{eqnarray}
P[\mathcal{W}_{1}(\mu_{n}^{N},\mu_{n}) >\varepsilon] &\leq& P[\mathcal{W}_{1}(\eta_{n}^{N},\mu_{n})> \gamma \varepsilon] +\sum_{i=0}^{n-1} P[\mathcal{W}_{1}(\eta_{i}^{N},\mu_{i}) \ge \frac{\gamma\varepsilon}{\delta M} (\frac{1-\gamma}{\chi_{2}})^{n-i}]\quad\quad\nonumber\\
&=& P[\mathcal{W}_{1}(\eta_{n}^{N},\mu_{n})> \gamma \varepsilon] +\sum_{i=1}^{n} P[\mathcal{W}_{1}(\eta_{n-i}^{N},\mu_{n-i}) \ge \frac{\gamma\varepsilon}{\delta M} \vartheta^{i}]\nonumber\\
&\equiv& T_1 + T_2,\label{unifexp1}
\end{eqnarray}
where $\vartheta = \frac{1-\gamma}{\chi_2}$, which, in view of our assumption on $\delta$, is strictly larger than $1$.
Let $i^{\e} = \max \{i \ge 0: \frac{\e \gamma}{\delta M}\vartheta^i < \frac{C_0}{2}\}$.
Then
$$
T_2 = \sum_{i=1}^{i^{\e}} P[\mathcal{W}_{1}(\eta_{n-i}^{N},\mu_{n-i}) \ge \frac{\gamma\varepsilon}{\delta M} \vartheta^{i}]
+ \sum_{i=i^{\e}+ 1}^{n} P[\mathcal{W}_{1}(\eta_{n-i}^{N},\mu_{n-i}) \ge \frac{\gamma\varepsilon}{\delta M} \vartheta^{i}].
$$
Note that since $\vartheta > 1$ and $t \mapsto \frac{\mathcal{C}^0(t)}{t^2}$ is non-increasing, $N \ge L_3 \frac{\log \mathcal{C}^0_{m_{\gamma}(\e)}}{m^2_{\gamma}(\e)}$ implies
$N \ge L_3 \frac{\log \mathcal{C}^0_{m_{\gamma}(\e \vartheta^i)}}{m^2_{\gamma}(\e \vartheta^i)}$ for all $i \ge 0$ where $m_{\gamma}$ is as introduced in Remark \ref{rem:vsconstant}.
Therefore from Lemma \ref{cor:cor621}(4), for all such $N$
$$
\sum_{i=1}^{i^{\e}} P[\mathcal{W}_{1}(\eta_{n-i}^{N},\mu_{n-i}) \ge \frac{\gamma\varepsilon}{\delta M} \vartheta^{i}] \le
\sum_{i=1}^{i^{\e}} \exp\{-L_2 N m^2_{\gamma}(\e \vartheta^i)\}.$$
Also, from Lemma \ref{cor:cor621}(3), for all $N \ge N_1$,
$$
\sum_{i=i^{\e}+ 1}^{n} P[\mathcal{W}_{1}(\eta_{n-i}^{N},\mu_{n-i}) \ge \frac{\gamma\varepsilon}{\delta M} \vartheta^{i}]
\le 
\sum_{i=i^{\e}+ 1}^{n} \exp\{-L_1 N m_{\gamma}(\e \vartheta^i)\}.$$
Combining these estimates and letting $N_2 = \max\{N_1, L_3\}$ and
$\tilde \vs_1(t) = \max\{1, \frac{\log \mathcal{C}^0_{m_{\gamma}(t)}}{m^2_{\gamma}(t)}\}$, we have for all $N \ge N_2 \tilde \vs_1(\e)$
$$
T_2 \le 2 \sum_{i=1}^{\infty} \exp \left \{ -L_4 N (\e^2\wedge \e)\vartheta^i\right \},$$
where
$L_4 = \min \{ L_2 \frac{\gamma^2}{ M^2}, L_1 \frac{\gamma}{ M}\}$.
Let $k_0 \in \mathbb{N}$ be such that for all $k \ge k_0$, $\vartheta^k \ge k$.  Then
$$
T_2 \le 2 k_0 \exp \left \{ -L_4 N (\e^2\wedge \e)\right \} + 2 \frac{\exp \left \{ -L_4 N (\e^2\wedge \e)\right \}}{1 - \exp \left \{ -L_4 N (\e^2\wedge \e)\right \}}.$$
Noting that $1 - \exp \left \{ -L_4 N (\e^2\wedge \e)\right \} \ge 1/2$ whenever $N \ge \frac{\log 2}{L_4}(\frac{1}{\e^2} \vee \frac{1}{\e})$, we see that
with
$\vs_1^*(t) = \max\{1, \frac{\log \mathcal{C}^0_{m_{\gamma}(t)}}{m^2_{\gamma}(t)}, \frac{1}{t^2}, \frac{1}{t}\}$
and $N_3 = \max\{N_1, L_3, \frac{\log 2}{L_4}\}$
\begin{equation}\label{eq:eq529b}
T_2 \le 2(k_0 + 2)\exp \left \{ -L_4 N (\e^2\wedge \e)\right \}\; \mbox{ for all } N \ge N_3 \vs_1^*(\e).
\end{equation}
Also from Lemma \ref{cor:cor621},
for all $N \ge N_3 \max \{ 1, \frac{\log \mathcal{C}^0_{\gamma \e}}{\gamma^2 \e^2}\}$
\begin{equation}\label{eq:eq529}
T_1 \le \exp \left \{ -L_5 N (\e^2\wedge \e)\right \}, \end{equation}
where $L_5 = \min \{\gamma^2 L_2, \gamma L_1\}$.
Using \eqref{eq:eq529b} and \eqref{eq:eq529} in \eqref{unifexp1} we now get the desired result with $a_1 = 2(k_0 +2) +1$, 
$a_2 = \min\{L_4, L_5\}$, $N_0 = N_3$ and $\vs_1(t) = \max\{\vs_1^*(t), \frac{\log \mathcal{C}^0_{\gamma t}}{\gamma^2 t^2}\}$.
\qed

\setcounter{section}{0}
 \setcounter{Theorem}{0}
 \setcounter{equation}{0}
 \renewcommand{\theequation}{\thesection.\arabic{equation}}
 \renewcommand {\theequation}{\arabic{section}.\arabic{equation}}

  \appendix
  \section*{Appendix}
  \renewcommand{\thesection}{A} 


The first part of the following lemma is an immediate consequence of Ascoli-Arzela theorem where as the second follows from Lemma 5 in \cite{como2009scaling}.
%
%
%
\begin{lemma}\label{app1}
(a) For a compact set $K$ in $\mathbb{R}^{d}$ let $\mathcal{F}_{a,b}(K)$ be the space of functions $f: K \to \mathbb{R}$ such that $\sup_{x\in K}|f(x)|\leq a$  and $|f(x) - f(y)|\leq b|x-y|$ for all $x,y \in K$. Then for any $\epsilon > 0$ there is a finite subset $\mathcal{F}_{a,b}^{\epsilon}(K)$ of $\mathcal{F}_{a,b}(K)$ such that for any signed measure $\mu$
$$\sup_{f \in \mathcal{F}_{a,b}(K)} |\langle f,\mu\rangle | \leq \max_{g \in \mathcal{F}_{a,b}^{\epsilon}(K)} |\langle g,\mu \rangle| + \epsilon |\mu|_{TV}.$$
(b) If $K = [-R, R]^d$ for some $R > 0$, then $\mathcal{F}_{R,1}^{\epsilon}(K)$ can be chosen such that 
$$|\mathcal{F}^{\epsilon}_{R,1}(K)| \leq \max\left \{\frac{2(2\sqrt{d}+1)}{3}  \frac{R}{\epsilon} 3^{[\frac{2R}{\epsilon}(\sqrt{d}+1)]^{d}}, 1 \right\}.$$
\end{lemma}
The next lemma is straightforward.
\begin{lemma}\label{app2}
Let $P:\mathbb{R}^{d}\times \mathcal{B}(\mathbb{R}^d) \to [0,1]$ be a transition probability kernel. Fix $N\geq 1$ and let $y_{1},y_{2},...,y_{N} \in \mathbb{R}^{d}$. Let $X_{1},X_{2},...,X_{N}$ be independent random variables such that $\mathcal{L}(X_{i})=\delta_{y_{i}}P.$ Let $f\in BM(\mathbb{R}^{d})$ and let $m_{0}^{N}= \frac{1}{N}\sum_{i=1}^{N}\delta_{y_{i}}$, $m_{1}^{N}= \frac{1}{N}\sum_{i=1}^{N}\delta_{X_{i}}$. Then 
$$E|\langle f, m_{1}^{N}- m_{0}^{N}P  \rangle| \leq \frac{2||f||_{\infty}}{\sqrt{N}}.$$
\end{lemma}
The following is a discrete version of Gronwall's lemma.
\begin{lemma}\label{app4}
Let $\{ a_{i}\}_{i=0}^{\infty},\{ b_{i}\}_{i=0}^{\infty},\{ c_{i}\}_{i=0}^{\infty}$ be  non-negative sequences. Suppose that
$$  a_{n}  \leq b_{n} + \sum_{k=0}^{n-1} c_{k} a_{k} \;   \mbox{ for all } n \ge 0.$$
Then 
$$ a_{n}  \leq b_{n} + \sum_{k=0}^{n-1}\left [ c_{k} b_{k}\left(\prod_{j=k+1}^{n-1}(1+ c_{j})\right) \right ]\;   \mbox{ for all } n \ge 0.$$
\end{lemma}


\skp

{\sc

\bigskip\noi
Amarjit Budhiraja\\
Department of Statistics and Operations Research\\
University of North Carolina\\
Chapel Hill, NC 27599, USA\\
email: budhiraj@email.unc.edu

\skp

\noi
Abhishek Pal Majumder\\
Department of Statistics and Operations Research\\
University of North Carolina\\
Chapel Hill, NC 27599, USA\\
email: palmajum@email.unc.edu

}

\end{document}